\documentclass[11pt]{amsart}
\usepackage{a4,amsmath,amssymb,amscd,verbatim}

\newtheorem{prop}[equation]{Proposition}

\newtheorem{thm}[equation]{Theorem}
\newtheorem{cor}[equation]{Corollary}
\newtheorem{lem}[equation]{Lemma}

\theoremstyle{definition}

\newtheorem{exa}[equation]{Example}

\numberwithin{equation}{section}

\newcommand{\spandsp}{\mbox{$\qquad\text{and}\qquad$}}

\newcommand{\sands}{\mbox{$\quad\text{and}\quad$}}

\newcommand{\tstts}[1]{\mbox{$\hspace{0.5em}\text{#1}\hspace{0.5em}$}}

\newcommand{\letbe}{\mathbin{:\!\raisebox{-.32pt}{=}}}

\newcommand{\hpm}{\mbox{$\hspace{.1em}$}}
\newcommand{\npm}{\mbox{$\hspace{-.1em}$}}

\newcommand{\ahss}{Atiyah-Hirzebruch spectral sequence}

\newcommand{\Sq}[1]{\mbox{${\it Sq\/}^{#1}$}}

\newcommand{\CP}{\text{\it CP\hspace{2pt}}}
\newcommand{\CPO}{\mbox{${\it CP\/}^1$}}
\newcommand{\CPT}{\mbox{${\it CP\/}^2$}}
\newcommand{\CPTp}{\mbox{${\it CP\/}^2_{\!+}$}}
\newcommand{\CPN}{\mbox{${\it CP\/}^n$}}
\newcommand{\CPNp}{\mbox{${\it CP\/}^n_{\!+}$}}
\newcommand{\CPI}{\mbox{${\it CP\/}^\infty$}}
\newcommand{\CPIp}{\mbox{${\it CP\/}^\infty_{\!+}$}}

\newcommand{\quare}{\scriptscriptstyle\square}
\newcommand{\tquare}{\mbox{$t_{\quare}$}}

\newcommand{\llongrightarrow}{\relbar\joinrel\longrightarrow}

\newcommand{\llongleftarrow}{\longleftarrow\joinrel\relbar}

\newcommand{\osu}{\mbox{$\varOmega_*^{U}$}}

\newcommand{\KO}{\mbox{\it KO\/}}

\newcommand{\BO}{\mbox{\it BO}}
\newcommand{\BSO}{\mbox{\it BSO}}
\newcommand{\BU}{\mbox{\it BU}}
\newcommand{\SO}{\mbox{\it SO}}

\newcommand{\bD}{\mathbb{D}}
\newcommand{\bC}{\mathbb{C}}
\newcommand{\bF}{\mathbb{F}}

\newcommand{\bQ}{\mathbb{Q}}
\newcommand{\bR}{\mathbb{R}}
\newcommand{\bZ}{\mathbb{Z}}

\newcommand{\brpk}{\mbox{$\bR_{\scriptscriptstyle\geqslant}^k$}}

\newcommand{\dlo}{\mbox{$d_{{\scriptscriptstyle\leq} 1}$}}
\newcommand{\dlj}{\mbox{$d_{{\scriptscriptstyle\leq} j}$}}
\newcommand{\dlk}{\mbox{$d_{{\scriptscriptstyle\leq} k}$}}
\newcommand{\velo}{\mbox{$v^E_{{\scriptscriptstyle\leq} 1}$}}
\newcommand{\velj}{\mbox{$v^E_{{\scriptscriptstyle\leq} j}$}}
\newcommand{\velk}{\mbox{$v^E_{{\scriptscriptstyle\leq} k}$}}

\newcommand{\srf}[1]{\mbox{${\text{\it SR}^F}$}}

\newcommand{\daaja}{Davis and Januszkiewicz}
\newcommand{\baabe}{Bahri and Benderskey}

\begin{document}
\bibliographystyle{plain}
\title[Decompositions and K--Theory of Bott Towers]
{Homotopy Decompositions and K--Theory\\
of Bott Towers}
\author{Yusuf Civan}
\address{Department of Mathematics, Faculty of Arts and Sciences, 
Suleyman Demirel University, 32260 Isparta, Turkey}
\email{ycivan@fef.sdu.edu.tr}
\author{Nigel Ray}
\address{Department of Mathematic, University of Manchester, Oxford
Road, Manchester M13~9PL, England}
\email{nige@ma.man.ac.uk}

\keywords {Bott towers, {\it K}-theory, stably complex structures,
Thom complexes, toric manifolds}


\begin{abstract}

We describe Bott towers as sequences of toric manifolds $M^k$, and
identify the omniorientations which correspond to their original
construction as toric varieties. We show that the suspension of $M^k$
is homotopy equivalent to a wedge of Thom complexes, and display its
complex $K$-theory as an algebra over the coefficient ring. We extend
the results to $KO$-theory for several families of examples, and
compute the effects of the realification homomorphism; these
calculations breathe geometric life into Bahri and Bendersky's
analysis of the Adams Spectral Sequence \cite{babe:ktt}. By way of
application we investigate stably complex structures on $M^k$,
identifying those which arise from omniorientations and those which
are almost complex. We conclude with observations on the r\^ole of
Bott towers in complex cobordism theory.

\end{abstract}

\maketitle

%
%
%
%
%
%
%
%
%

\section{Introduction}\label{intro}

In their 1950s study of loops on symmetric spaces, Bott and Samelson
\cite{bosa:atm} introduced a remarkably rich and versatile family of
smooth manifolds. Various special cases were treated in different
contexts during the following three decades, until Grossberg and Karshon
\cite{grka:btc} offered a description as complex algebraic varieties in
1994. They referred to their constructions as {\it Bott towers}, and
addressed issues of representation theory and symplectic geometry. Our
purpose here is to offer the alternative viewpoint of algebraic
topology. We consider Bott towers $(M^k:k\leq n)$ of height $n$, and
discuss homotopy decompositions of the suspensions $\varSigma M^k$;
these provides further evidence that the spaces of complex geometry are
often stably homotopy equivalent to wedges of Thom complexes, as we have
argued elsewhere \cite{grra:scb}. We investigate the real and complex
$K$-theory of the $M^k$, casting geometric light on recent calculations
of Bahri and Bendersky \cite{babe:ktt} which were originally conducted
in the algebraic underworld of the Adams Spectral Sequence.

Given a commutative ring spectrum $E$, we denote the reduced and
unreduced cohomology algebras of any space $X$ by $E^*(X)$ and
$E^*(X_+)$ respectively. So $E^*(S^n)$ is a free module over the
coefficient ring $E_*$ on a single $n$-dimensional generator $s^E_n$,
defined by the unit of $E$. In particular, we use this notation for the
integral Eilenberg-Mac Lane spectrum $H$ and the complex $K$-theory
spectrum $K$. Real $K$-theory requires the most detailed calculations,
so we abbreviate $s^{KO}_n$ to $s_n$ whenever possible. We require
multiplicative maps $f\colon E\rightarrow F$ of ring spectra to preserve
the units, so that $f(s^E_n)=s^F_n$ for all $n$; complexification
$c\colon\KO\rightarrow K$ is an important example. We adopt similar
conventions for Thom classes $t^E$, which also play a major
r\^ole. Given an $E$-orientable $n$--dimensional vector bundle $\gamma$,
we insist that $t^E$ should lie in $E^n(T(\gamma))$, and restrict to
$s_n^E$ on the fibre. Alternative choices of dimension are, of course,
available for periodic spectra such as $K$ and $\KO$, but we believe
that our chosen convention leads to the least confusion.

With the single exception of $\KO$, the spectra we use are complex
oriented by an appropriate choice of first Chern class $v^E$ in
$E^2(\CPI)$; by definition, $v^E$ restricts to $s^E_2$ on $\CPO$. We
also insist that $E_*$ be concentrated in even degrees.

The contents of our sections are as follows.

In Section \ref{twgeco} we establish our notation, and recall
well-known computations for the $E$-cohomology of certain sphere
bundles $Y$ over complexes with cells in even dimensions. We record a
homotopy decomposition of $\varSigma Y$, and apply the results to
$K$-theory and integral cohomology. We introduce Bott towers as
iterated sphere bundles in Section \ref{boto}, and apply the previous
section to describe their $E$-cohomology algebras, and splittings of
their suspensions. We also consider their stable tangent bundles, and
introduce a cofiber sequence relating pairs of towers. Bott towers
masquerade as toric varieties, and we discuss their associated
properties in Section \ref{tost}; we adapt the viewpoint of Grossberg
and Karshon, and pay particular attention to the corresponding complex
structures. Our calculations with $\KO$-theory begin in Section
\ref{rekstontw}, where we focus on dimensions $2$ and $4$. We obtain
complete descriptions of the $\KO_*$-algebra structure in all cases.
These results provide a springboard for our most comprehensive
calculations, which occupy Section \ref{rekboto}; we consider all
dimensions, but specialise to two particular families of cases. Again,
we obtain complete information about $\KO_*$-algebra structures, but
find that certain products are particularly complicated to describe
explicitly. We relate our results to the pioneering work of Bahri and
Bendersky. Finally, in Section \ref{stcost}, we apply these
calculations to the enumeration of a collection of stably complex
structures, which arise from our study of Bott towers as toric
manifolds. Such structures are of key importance to understanding
their r\^ole in complex cobordism theory.

The idea of studying Bott towers in this context first emerged during
discussions with Victor Buchstaber, made possible by Aeroflot's
abandonment of flights out of Manchester in 1996. The second author
announced most of the results at the Conference on Algebraic Topology
in Gdansk, Poland, during June 2001, where Taras Panov and his
colleagues offered many helpful suggestions as we strolled the Baltic
beaches. We apologise to them all for our protracted attempts to
produce a final document, and give thanks to Adrian Dobson for
identifying several errors in various intermediate versions.

%
%
%
%
%
%
%
%
%

\section{2--Generated Complexes}\label{twgeco}

It is convenient to work with connected CW-complexes $X$ whose
integral cohomology ring $H^*(X;\bZ)$ is generated by a linearly
independent set of $2$--dimensional elements $x_1$, \dots, $x_m$. We
describe such an $X$ as being {\it 2--generated}, and note that
$H^2(X;\bZ)$ is isomorphic to the integral lattice $\bZ^m$; we refer
to the elements $x_j$ as the {\it 2--generators\/} of $X$, and to $m$
as its {\it 2--rank}. We follow combinatorial convention by
abbreviating the set $\{1,\dots,m\}$ to $[m]$, and denote the product
$\prod_Rx_j$ by $x_R$ for any subset $R\subseteq[m]$. The first Chern
class $v^H$ defines a canonical isomorphism between the multiplicative
group of complex line bundles over $X$ and $H^2(X;\bZ)$, and so
determines line bundles $\gamma_j$ such that $v^H(\gamma_j)=x_j$, for
$1\leq j\leq m$. In general, it assigns the $m$-tuple
$(a(1),\dots,a(m))$ to the tensor product
\begin{equation}\label{tenbun}
\gamma_1^{a(1)}\otimes\dots\otimes\gamma_m^{a(m)}.
\end{equation}
By definition, $X$ lies in the category of CW-complexes whose cells
are even dimensional. Various observations of Hoggar \cite{ho:ktg}
therefore apply to the abelian group structure of $\KO^*(X)$, and are
relevant to parts of Sections \ref{rekstontw} and \ref{rekboto}.

Given any of our complex oriented ring spectra $E$, the Chern classes
$v^E(\gamma_j)=v^E_j$ lie in $E^2(X)$ for all $1\leq j\leq m$. The
corresponding \ahss\ collapses for dimensional reasons, and identifies
$E^*(X)$ as a free $E_*$-module, spanned by the monomials $v^E_R$; in
other words, it is generated by $v^E_1$, \dots, $v^E_m$ as an
$E_*$-algebra. An important, if atypical, example is provided by
$\CPN$. Then $v=v^E(\zeta(n))$ is the first Chern class of the Hopf
line bundle $\zeta(n)$, and the canonical isomorphism
\begin{equation}\label{ecpn}
E^*(\CPNp)\;\cong\;E_*[[v]]/(v^{n+1})
\end{equation}
confirms that $\CPN$ has the single $2$--generator $v$. In order to
emphasise that we are working over $\CPN$, we sometimes denote $v$ by
$v(n)$; thus $v(1)$ and $s^E_2$ are interchangeable. In the cases
$E=H$ and $K$, we write $v$ as $x$ and $u$ respectively.

The following results are well-known, and are usually obtained by
applying standard methods of Borel and Hirzebruch \cite{bohi:cch}. Our
immediate interests, however, are homotopy theoretic, and involve the
stable triviality of certain cofibre sequences of $2$--generated
complexes and associated Thom spaces. We therefore take the
opportunity to establish our notation by outlining proofs in this
alternative language.

We assume that $X$ is 2--generated, and write $\gamma$ for the line
bundle \eqref{tenbun}. We let $Y$ denote the total space
$S(\bR\oplus\gamma)$ of the 2--sphere bundle obtained from $\gamma$
by the addition of a trivial real line bundle, and write $p$ for the
projection onto $X$. Whenever $X$ is a smooth manifold, we may
assume that $Y$ is also.
\begin{lem}\label{twogen}
The $E_*$-algebra $E^*(Y_+)$ is a free module over $E^*(X_+)$ on
generators $1$ and $v^E_{m+1}$, which have dimensions $0$ and $2$
respectively; the multiplicative structure is determined by the
single relation
\begin{equation}\label{multreln}
(v^E_{m+1})^2=v^E(\gamma)v^E_{m+1},
\end{equation}
and $v^E_{m+1}$ restricts to $s^E_2$ on the fibre $S^2\subset Y$.
\end{lem}
\begin{proof}
The sphere bundle $S(\bR\oplus\gamma)$ admits a section $r$, given
by $+1$ in the summand $\bR$, and the quotient of the total space by
the image of $r$ is canonically homeomorphic to the Thom complex
$T(\gamma)$ \cite{st:nct}. In the resulting cofibre sequence
\begin{equation}\label{xycofseq}
X\stackrel{r}{\longrightarrow}Y\stackrel{q}{\longrightarrow}T(\gamma),
\end{equation}
the quotient map $q$ identifies the fibres $S^2\subset Y$ and
$S^2\subset T(\gamma)$, and $r$ has left inverse $p$. The standard
coaction of $X$ on $T(\gamma)$ interacts with the diagonal on $Y$ by
the commutative square
\begin{equation}\label{multstr}
\begin{CD}
Y@>>q> T(\gamma)\\ @V\delta VV@VV\delta V\\ Y\times
Y@>(p,\hspace{.3pt}q)>>X_+\wedge T(\gamma)
\end{CD}.
\end{equation}

The $E$-cohomology sequence induced by \eqref{xycofseq} is split by
$p^*$, and is therefore short exact. The Chern class $v^E$ induces a
canonical Thom class $t^E\in E^2(T(\gamma))$, and so determines a
Thom isomorphism $E^{*-2}(X_+)\cong E^*(T(\gamma))$, which
identifies $E^*(Y_+)$ as the free $E^*(X_+)$--module on generators
$1$ and $v^E_{m+1}=q^*t^E$. The diagram \eqref{multstr} confirms
that products of the form $p^*(x)v^E_{m+1}$ may be written as
$q^*(xt^E)$ for any $x\in E^*(X)$; so the action of $v^E_{m+1}$ is
by multiplication in $E^*(Y)$. Since $\delta^*(v^E(\gamma)\otimes
t^E)=(t^E)^2$, the formula for $(v^E_{m+1})^2$ follows.
\end{proof}

An obvious consequence of Lemma \ref{twogen} is that $Y$ is also
2--generated, and has 2--rank $m+1$. The Chern class $v^E(\gamma)$
may be expanded in terms of the $E_*$-basis $v^E_1$, \dots, $v^E_m$,
using the associated formal group law $F^E$ and its $n$-series
$[n]^E$. We obtain
\begin{equation}\label{fgp}
v^E(\gamma)=F^E([a(1)]^E,\dots,[a(m)]^E)
\end{equation}
in $E^2(x)$.

The universal example of Lemma \ref{twogen} is given by $X=\CPI$ and
$\gamma=\zeta$; it follows that $T(\gamma)$ is also homeomorphic to
$\CPI$, and that $Y$ is homotopy equivalent to $\CPI\vee\CPI$. Then
$E^*(Y_+)$ is free over $E_*[[v]]$ on generators $1$ and $v'$, with
$(v')^2=vv'$. The general case may be deduced from this example by
pulling back along the classifying map for $\gamma$.  Of course, we
may restrict the universal example to any skeleton $X=\CPN$, in which
case $T(\gamma)$ is ${\it CP\/}^{n+1}$.

There is a second section $\widetilde{r}\colon X\rightarrow Y$,
defined by $-1\in\bR$. The resulting composition
$q\cdot\widetilde{r}\colon X\rightarrow T(\gamma)$ reduces to the
inclusion of the zero-section, giving
$\widetilde{r}^{\,*}t^E=v^E(\gamma$).

The usual approach to Lemma \ref{twogen} proceeds by identifying
$S(\bR\oplus\gamma)$ with its {\it projective form} ${\it
CP}(\bC\oplus\gamma)$. The corresponding canonical line bundle has
first Chern class $v^E_{m+1}$, and is isomorphic to $\gamma_{m+1}$;
it restricts to the Hopf bundle $\zeta(1)$ over the fibre $\CPO$. So
$\gamma_{m+1}$ is a summand of the pullback $\bC\oplus\gamma$ over
$Y$, and has orthogonal complement
$\overline{\gamma}_{m+1}\otimes\gamma$ with respect to the standard
inner product. The associated splitting
\begin{equation}\label{gamgamdav} \bC\oplus
p^*\gamma\;\cong\;
\gamma_{m+1}\oplus\left(\overline{\gamma}_{m+1}\otimes
p^*\gamma\right)
\end{equation}
gives rise to the relation \eqref{multreln}, and will be useful
in Section \ref{stcost}.

The cofibre sequence \eqref{xycofseq} also leads to the familiar
relationship between the homotopy types of $X$ and $Y$.
\begin{prop}\label{xysusp}
There is a homotopy equivalence
\[
h\colon\varSigma Y\longrightarrow\varSigma X\vee\varSigma T(\gamma)
\]
of suspensions.
\end{prop}
\begin{proof}
We define $h$ as the sum $\varSigma p+\varSigma q$, and construct a
homotopy inverse $\varSigma X\vee\varSigma
T(\gamma)\rightarrow\varSigma Y$ by forming the wedge of $\varSigma
r$ with the map $l\colon\varSigma T(\gamma)\rightarrow\varSigma Y$
which collapses the standard copy of $X$ in $T(\bR\oplus\gamma)$.
\end{proof}
The equivalence $h$ induces an isomorphism in $E$-cohomology, which
realises the module structures of Lemma \ref{twogen} by splitting
$E^*(Y_+)$ as $E^*(X_+)\oplus(v^E_{m+1})$. In the universal example,
$h$ is a self equivalence of $\varSigma\CPI\vee\varSigma\CPI$ and 
desuspends.

We shall need an extension of Lemma \ref{twogen}, in the situation
when $X$ itself is the total space of a bundle $\theta$ over $S^2$,
with fibre $X'$. We write $\gamma'$ for the pullback of $\gamma$ to
$X'$, and $Y'$ for the total space $S(\bR\oplus\gamma')$; thus $Y'$
is also the fibre of the projection $Y\rightarrow S^2$.
\begin{prop}\label{overstwo}
With the data above, there is a homotopy commutative ladder of
cofibre sequences
\begin{equation}\label{cladd}
\begin{CD}
T(\gamma')@>i>>T(\gamma)@>f>>\varSigma^2T(\gamma')\\
@Aq'AA @AqAA @AA\varSigma^2q'A\\
Y'@>>i>Y@>>f>\varSigma^2Y'_+
\end{CD}\;,
\end{equation}
where the maps $i$ are induced by inclusion of the fibre, and the maps
$f$ are quotients.
\end{prop}
\begin{proof}
We may construct $X$ from two copies of $\bD^2\times X'$ by
identifying them along their boundaries $S^1\times X'$ via the
characteristic function of $\theta$. Then $X/i(X')$ is homeomorphic
to $\varSigma^2X'_+$. The same argument applies to $Y/i(Y')$,
yielding cofibre sequences
\begin{equation}\label{xpxcofseq}
X'\stackrel{i}{\longrightarrow}X
\stackrel{f}{\longrightarrow}\varSigma^2X'_+\spandsp
Y'\stackrel{i}{\longrightarrow}Y
\stackrel{f}{\longrightarrow}\varSigma^2Y'_+.
\end{equation}
The sections $r'\colon X'\rightarrow Y'$ and $r\colon X\rightarrow Y$
are compatible with the inclusions $i$, and the ladder follows by
taking quotient maps $q'$ and $q$.
\end{proof}

The naturality of the ladder \eqref{cladd} leads to a commutative square
\begin{equation}\label{squdiag}
\begin{CD}
T(\gamma)@>>f>\varSigma^2T(\gamma')\\
@V\delta VV@VV\epsilon V\\
X\wedge T(\gamma)@>f\wedge 1>>\varSigma^2(X'_+)\wedge T(\gamma)
\end{CD},
\end{equation}
where $\epsilon$ is the Thom complexification of the bundle map
obtained by pulling $\bR^2\times\gamma$ back along the restricted
diagonal $X'\rightarrow X'\times X$. Alternatively, the square may be
considered as the quotient of the reduced diagonal
$T(\gamma)\rightarrow X\wedge T(\gamma)$ by its restriction
$T(\gamma')\rightarrow X'\wedge T(\gamma)$.

The first sequence of \eqref{xpxcofseq} induces the Wang long exact
sequence of $\theta$ in $E$-cohomology, for any multiplicative
spectrum $E$. Standard homotopy theoretic arguments \cite{wh:eht}
show that the connecting map $\varSigma^2X'_+\rightarrow\varSigma
X'$ is induced from the characteristic map $S^1\times X'\rightarrow
X'$ by suspension.

We shall apply these facts in the particular cases $E=H$ and $K$,
denoting the elements $v^E_{m+1}$ by $x_{m+1}$ and $g_{m+1}$
respectively. We write the coefficients of complex $K$-theory as the
ring of Laurent series
\[
K_*\letbe\bZ[z,z^{-1}],
\]
where $z$ lies in $K_2$ and is represented by the virtual Hopf line
bundle over $S^2$. So $zg_j$ is represented by the virtual
bundle $\gamma_j-\bC$ in $K^0(X)$, for $1\leq j\leq m$. Complex
conjugation acts on $K_*$ by $\overline{z}=-z$, and on the algebra
generators by
\begin{equation}\label{conjongj}
\overline{g}_j\;=\;\overline{\gamma}_jg_j\;=\;g_j/(1+zg_j)\;=\;
\;\sum_{i=0}^{\infty}(-z)^ig_j^{i+1};
\end{equation}
the Chern character embeds $K^*(X)$ in the ring
$H^*(X;\bQ[z,z^{-1}])$ by $ch(g_j)=z^{-1}(e^{zx_j}-1)$, for $1\leq
j\leq m$.

The cases $H$ and $K$ correspond to the additive and multiplicative
formal group laws respectively. The Chern classes \eqref{fgp} are
given by
\begin{equation}\label{eehak}
\begin{split}
v^H(\gamma)=a(1)&x_1+\dots+a(m)x_m\sands\\
&v^K(\gamma)=z^{-1}\Big(\prod_{j\leq m}(1+zg_j)^{a(j)}-1\Big),
\end{split}
\end{equation}
and are compatible under the action of the Chern character.

%
%
%
%
%
%
%
%
%

\section{Bott Towers}\label{boto}

In this section we consider the algebraic topology of Bott towers,
extending our work \cite{ra:ocb} on bounded flag manifolds; our
methods complement the more geometric approach of \cite{ci:phd}. We
give an inductive construction as a family of 2--generated smooth
oriented manifolds $M^k$, and describe their cohomology rings for any
of our complex oriented ring spectra $E$. We obtain an elementary
decomposition of their suspensions $\varSigma M^k$ into a wedge of
Thom complexes, and consider two natural complex structures on their
stable tangent bundles.

Given any integer $k\geq 1$, we assume that a {\it $(k-1)$th
stage\/} $M^{k-1}$ has been constructed as a smooth oriented
$2(k-1)$--dimensional manifold with 2--generators $v^E_j$, and line
bundles $\gamma_j$ such that $v^E(\gamma_j)=v^E_j$, for $1\leq j\leq
k-1$. Using the notation of \eqref{tenbun}, we write
$\gamma(a_{k-1})$ for the complex line bundle
\[
\gamma_1^{a(1,k)}\otimes\dots\otimes\gamma_{k-1}^{a(k-1,k)}
\]
associated to the $(k-1)$-tuple $a_{k-1}=(a(1,k),\dots,a(k-1,k))$ in
$\bZ^{k-1}$. Fixing $a_{k-1}$, we refer to $\gamma(a_{k-1})$ as the
$k${\it th bundle\/} of the construction, and define $M^k$ to be the
total space of the smooth 2--sphere bundle of
$\bR\oplus\gamma(a_{k-1})$, oriented by the outward pointing normal
and the complex structure on $\gamma(a_{k-1})$. By Lemma
\ref{twogen}, we deduce that $M^k$ has 2--generators $v^E_j$ for
$1\leq j\leq k$, where $v^E_k$ is the pullback of the Thom class
$t^E_k\in E^2(T(\gamma(a_{k-1})))$ along the collapse map $q_k$.
Moreover, $t^E_k$ is the first Chern class of a canonical line
bundle $\lambda_{k-1}$ over $T(\gamma(a_{k-1}))$, so
$v^E_k=v^E(\gamma_k)$, where $\gamma_k$ is defined as
$q_k^*\lambda_{k-1}$. Henceforth, we abbreviate $T(\gamma(a_j))$ to
$T(a_j)$ for each $1\leq j\leq n$.

In order to get off the ground, it is convenient to write the
one-point space as $M^0$, so that the first bundle is trivial and
$x_0=0$. Then $M^1$ is a $2$--sphere, compatibly oriented with the
complex structure on $\CPO$, and $\gamma_1$ is the Hopf line bundle
$\zeta(1)$. The cohomology ring $E^*(S^2_+)$ is isomorphic to
$E_*[v]/(v^2)$, where $v=v^E(\zeta(1))=s^E_2$, and $S^2$ is
$2$--generated with $2$--rank $1$. Of course the second bundle
$\gamma(a_1)$ is isomorphic to $\zeta(1)^{a(1,2)}$ for some $1$-term
sequence $a_1=(a(1,2))$.

The construction is now complete, and the $k$th stage depends only
on the integral sequences $(a_1,\dots,a_{k-1})$, which contain
$k(k-1)/2$ integers $a(i,j)$, for $1\leq i<j\leq k-1$. It is
occasionally helpful to interpret $a_0$ as empty, and to write the
first bundle $\bC$ as $\gamma(a_0)$.

We refer to the sequence $(M^k:k\leq n)$ of oriented manifolds as a
{\it Bott tower\/} of {\it height} $n$ (which may be infinite); it is
determined by the {\it list} $a=(a_1,\dots,a_{n-1})$ of $n(n-1)/2$
integers. If we choose the projective form of $M^k$ at every stage, we
obtain a tower of nonsingular algebraic varieties, whose orientations
coincide with those decribed above. Every Bott tower involves
projections $p_k\colon M^k\rightarrow M^{k-1}$, sections $r_k$ and
$\widetilde{r}_k\colon M^{k-1}\rightarrow M^k$, and quotient maps
$q_k\colon M^k\rightarrow T(a_{k-1})$, for each $1\leq k\leq n$.

The cohomological structure of $M^k$ is given as follows.
\begin{prop}\label{cohomk}
For any complex oriented ring spectrum $E$, the $E_*$-algebra
$E^*(M^k_+)$ is isomorphic to $E_*[v^E_1,\dots,v^E_k]/I^E_k$, where
$I^E_k$ denotes the ideal
\[
\big((v^E_j)^2-v^E(\gamma)v^E_j:1\leq j\leq k\big);
\]
in particular, $E^{2r}(M^k_+)$ is the free $E_*$-module generated by
the monomials $v^E_R$, as $R\subseteq[k]$ ranges over the subsets of
cardinality $r$, and $E^*(M^k_+)$ has total rank $2^r$.
\end{prop}
\begin{proof}
The multiplicative structure follows from $k-1$ applications of
Lemma \ref{twogen}; the resulting relations imply the additive
structure immediately.
\end{proof}

In the cases $E=H$ and $K$, we denote the elements $v^E_j$ by $x_j$ in
$H^2(M^k_+;\bZ)$ and $g_j$ in $K^2(M^k_+)$ respectively, for $1\leq
j\leq k$. The ideals $I^H_k$ and $I^K_k$ are then described explicitly
by \eqref{eehak}. The structure of $H^*(M^k_+;\bZ)$ shows that the
Euler characteristic of $M^k$ is $2^k$, and is independent of $a$;
this may also be confirmed by straightforward geometric argument.

By way of example we consider the tower $(B_k:0\leq k)$, whose list
satisfies $a_k=(0,\dots,0,1)$ for all $k\geq 1$. We studied this
example in \cite{ra:ocb}, where we explained its significance for
complex cobordism theory. In later work \cite{bura:fml} we interpreted
the points of $B_k$ as complete flags $0<U_1<\dots<U_n<\bC^{k+1}$,
bounded below by the standard flag in the sense that the first $j$
standard basis vectors lie in $U_{j+1}$, for each $1\leq j\leq k$. The
resulting description of $B_k$ as a {\it bounded flag manifold\/}
corresponds to the projective form ${\it CP}(\bC\oplus\gamma_{k-1})$,
and displays $B_k$ as a toric variety.

We may now describe our homotopy theoretic decomposition of
$\varSigma M^k$.
\begin{prop}\label{btsplit}
Given any Bott tower $(M^k:k\leq n)$, there is a homotopy
equivalence
\[
h_k\colon \varSigma M^k\longrightarrow\varSigma S^2\vee\varSigma
T(a_1)\vee\dots\vee\varSigma T(a_{k-1}),
\]
for each $1\leq k\leq n$.
\end{prop}
\begin{proof}
It suffices to apply Proposition \ref{xysusp} $k-1$ times; $S^2$
appears as $T(a_0)$.
\end{proof}

With respect to Proposition \ref{cohomk}, the homotopy equivalence
$h_k$ induces the additive splitting
\[
E^*(M^k)\;\cong\;
\langle\velo\rangle\oplus\dots\oplus\langle\velk\rangle
\]
where $\langle\velj\rangle$ denotes the free $E_*$-submodule
generated by those monomials $v^E_R$ for which $R\subseteq[j]$ and
$j\in R$. By construction, $\langle\velj\rangle$ is the image of
$E^*(T(a_{j-1}))$ under the injection $p_k^*\cdots p_{j+1}^*q_j^*$,
for each $1\leq j\leq k$; it is split by $l_{j-1}^*r_j^*\cdots
r_k^*$, where $l_{j-1}^*$ is induced by the map $\varSigma
T(a_{j-1})\rightarrow\varSigma M^{j-1}$ which collapses the standard
copy of $M^{j-1}$ in $T(\bR\oplus\gamma(a_{j-1}))$.

It is worth commenting on aspects of the case $k=2$, which is influenced
by the fact that the isomorphism class of the {\it SO}(3)-bundle
$\bR\oplus\zeta(1)^{a(1,2)}$ depends only on the parity of $a(1,2)$. So
there are diffeomorphisms $ M^2\rightarrow S^2\times S^2$ when
$a(1,2)=2b$ is even, and $M^2\rightarrow S(\bR\oplus\zeta(1))$ when
$a(1,2)=2b+1$ is odd. In $E$-cohomology, they induce isomorphisms
\begin{equation}\label{mtisos}
\begin{split}
E_*[v_1,v_2]/\left(v_1^2,v_2^2-2b
v_1v_2\right)&\;\cong\;
E_*[w_1,w_2]/(w_1^2,w_2^2)\sands\\E_*
[v_1,v_2]/\left(v_1^2,v_2^2-(2b+1)v_1v_2\right)&\;\cong\;
E_*[w_1,w_2]/\left(w_1^2,w_2^2-w_1w_2\right),
\end{split}
\end{equation}
(omitting the superscripts $E$), which are determined by the $2\times
2$ matrices of their actions on the column vector $(v_1,v_2)$. Such
matrices are exemplified by
$\left(\begin{smallmatrix}1&0\\b&1\end{smallmatrix}\right)$, for any
integer $b$.

We shall be particularly interested in the stable tangent bundle of
$M^k$ in Section \ref{stcost} below. As explained by Szczarba
\cite{sz:tbf}, there is an explicit isomorphism
\begin{equation}\label{stanmk}
\tau(M^k)\oplus\bR\;\cong\;
\bR\oplus\,\bigoplus_{j=1}^k\gamma(a_{j-1})
\end{equation}
of $\SO(2k+1)$--bundles, which determines a stably almost complex
structure $\tau'$ on $M^k$. Since \eqref{stanmk} extends over the
$3$--disk bundle of $\bR\oplus\gamma(a_{k-1})$, this structure
bounds. On the other hand, the projective form of $M^k$ is a
nonsingular complex algebraic variety, whose tangent bundle admits the
canonical complex structure described in Section \ref{tost}. The fact
that its stabilisation differs from \eqref{stanmk} is one of our
motivations for Section \ref{stcost}.

Given a Bott tower of height $n$, we turn our attention to the
projection $p_{n,k}\colon M^n\rightarrow M^k$, defined as the
composition $p_{k+1}\cdots p_n$ for some $k\geq 1$. This is also a
smooth bundle, whose fibre we wish to identify.
\begin{prop}\label{newfib}
The fibre of $p_{n,k}$ is the $(n-k)$th stage of a Bott tower
$((M')^j:j\leq n-k)$; it is determined by the list
$(a_1',\dots,a_{n-k-1}')$, where $a_j'$ is formed from $a_{j+k}$ by
deleting the first $k$ entries, for each $1\leq j\leq n-k-1$.
\end{prop}
\begin{proof}
When we restrict the bundle $\gamma(a_k)$ to a point $(M')^0$ in
$M^k$, we obtain the trivial bundle $\bC$, and $M^{k+1}$ pulls back
to the fibre $S^2$ of $p_{k+1}$; we label this fibre $(M')^1$. We
repeat the pullback procedure over $(M')^1$, and continue until we
reach $M^{n-1}$. We find that $\gamma_j$ restricts trivially to
$(M')^{n-k-1}$ for $1\leq j\leq k$, and to $\gamma_{j-k}'$ for
$k<j\leq n-1$. Thus $\gamma(a_{n-1})$ restricts to
$\gamma'(a_{n-k-1}')$, where $a_{n-k-1}'=
(a(k+1,n),\dots,a(n-1,n))$, and $M^n$ pulls back to
$S(\bR\oplus\gamma'(a_{n-k-1}'))$, which we label $(M')^{n-k}$. The
construction ensures that $(M')^{n-k}$ is the inverse image of
$(M')^0$ under $p_{n,k}$, and is therefore the required fibre.
\end{proof}
\begin{cor}\label{decomptwo}
For each $1<k<n$, there is a commutative ladder of cofibre sequences
\begin{equation}\label{ladder}
\begin{CD}
T(a_{k-2}')@>i>>T(a_{k-1})@>f>>\varSigma^2T(a_{k-2}')\\ @Aq_{k-1}'AA
@Aq_kAA @AA\varSigma^2q_{k-1}'A\\
(M')^{k-1}@>>i>M^k@>>f>\varSigma^2(M')^{k-1}_+
\end{CD}\quad.
\end{equation}
In $E$-cohomology, the homomorphisms induced by the upper sequence
satisfy $i^*t^E_k=(t')^E_{k-1}$, and
$f^*(\varSigma^2i^*w(t')^E_{k-1})=v^E_1wt^E_k$ for every $w\in
E^*(M^{k-1})$. In the lower sequence they satisfy
$i^*v^E_j=(v')^E_{j-1}$ for each $2\leq j\leq k$, with $i^*v^E_1=0$,
and $f^*(\varSigma^2i^*v^E_R)=v^E_1v^E_R$ for every
$R\subseteq\{2,\dots,k\}$, with $f^*s^E_2=v^E_1$.
\end{cor}
\begin{proof}
The ladder arises by combining Proposition \ref{newfib} with
Proposition \ref{overstwo}, where $X$ is $M^{k-1}$ and $Y$ is $M^k$.
Since the upper $i$ arises from a bundle map it satisfies
$i^*t^E_k=(t')_{k-1}^E$, yielding $i^*v^E_k=(v')_{k-1}^E$; the
corresponding result holds for $j<k$ by projection onto $M^j$, noting
that $(t')^E_0=0$. Pulling $s^E_2\otimes wt^E_k$ back around
\eqref{squdiag} confirms that
$f^*(\varSigma^2i^*w(t')^E_{k-1})=v^E_1wt^E_k$ in $E^*(T(a_{k-1})$,
and applying \eqref{multstr} leads to the formula for $f^*$ on
$E^*(\varSigma^2(M')^{k-1}_+)$.
\end{proof}

Since all the spaces on view in Corollary \ref{decomptwo} are
$2$-generated, the horizontal cofibre sequences are cohomologically
split. The formulae for $i^*$ and $f^*$ show that the splitting of
$E^*(M^k_+)$ take the form
\begin{equation}\label{split1}
\begin{split}
E_*[&v^E_1,\dots,v^E_k]/I^E_k\;\cong\;\\
&E_*[(v')_1^E,\dots,(v')_{k-1}^E]/(I')^E_{k-1})\oplus
v^E_1E_*[(v')_1^E,\dots,(v')_{k-1}^E]/(I')^E_{k-1},
\end{split}
\end{equation}
and subsumes the splitting of $E^*(T(a_{k-1}))$ as
\[
\langle\velk\rangle\cong\langle
(v')_{{\scriptscriptstyle\leq}k-1}^E\rangle \oplus v^E_1\langle
(v')_{{\scriptscriptstyle\leq}k-1}^E\rangle.
\]

%
%
%
%
%
%
%
%
%

\section{Toric Structures}\label{tost}

We now describe the stages of a Bott tower $(M^k:k\leq n)$ as toric
manifolds, in the sense of \daaja; we continue to assume that the
tower is determined by the list $a=(a_1,\dots,a_{n-1})$. We use the
language of \cite{bura:tst} to record the salient properties, and
discuss the relationship with Grossberg and Karshon's construction
\cite{grka:btc} of the $M^k$ as complex manifolds.

We write the $k$-dimensional torus as $T^k$ and denote a generic
point $t$ by $(t_1,\dots,t_k)$, where $t_i$ lies in the unit circle
$T\subset\bC$ for each $1\leq i\leq k$. So $T^k$ is naturally
embedded in $\bC^k$, on which it acts coordinatewise, by
multiplication; this is the {\it standard action\/}, whose quotient
is the nonnegative orthant $\brpk$. We study the standard action of
$T^{2k}$ on $(S^3)^k$, induced by embedding the latter in $\bC^{2k}$
as the subspace
\begin{equation}\label{sthreek}
\{(y_1,z_1,\dots,y_k,z_k):
y_i\overline{y}_i+z_i\overline{z}_i=1\tstts{for}1\leq i\leq k\}.
\end{equation}
When $k=1$ the quotient of this action is a curvilinear $1$--simplex, or
interval, $I=\{(r,s):r^2+s^2=1\}$ in
$\bR_{\scriptscriptstyle\geqslant}^2$, so for general $k$ it is a
curvilinear cube $I^k\subset\bR_{\scriptscriptstyle\geqslant}^{2k}$.

Given $a$, we define the $k$--dimensional subtorus $T^k(a)<T^{2k}$ to
consist of elements
\begin{equation}\label{torusta}
\begin{split}
\{(u_1,u_1,u_2,u_1^{-a(1,2)}u_2,\dots,u_k,u_1^{-a(1,k)}\dots
&u_{k-1}^{-a(k-1,k)}u_k):\\
&u_i\in T\tstts{for}1\leq i\leq k\},
\end{split}
\end{equation}
for each $k\leq n$. So $T^k(a)$ acts freely on $(S^3)^k$, and the
quotient space $Q_k$ is a smooth $2k$--dimensional manifold. Moreover,
the $k$--torus $T^{2k}/T^k(a)$ acts on $Q_k$, and has quotient $I^k$;
with respect to this action, $Q_k$ is a toric manifold. We abbreviate
$T^{2k}/T^k(a)$ to $T^k_a$ whenever it acts on $Q_k$ in this fashion.

\begin{prop}\label{qisproj}
Given any $1\leq k\leq n$, there is an orientation preserving
diffeomorphism $\phi_k\colon Q_k\rightarrow M^k$; it pulls $\gamma_j$
back to the line bundle 
\[
(S^3)^k\times_{T^k(a)}\bC\longrightarrow Q_k
\]
for each $1\leq j\leq k$, where $T^k(a)$ acts on $\bC$ by
$w\mapsto u_j^{-1}w$.
\end{prop}
\begin{proof}
We proceed by induction on $k$, noting that $\phi_1$ is defined by
factoring out the action of $T^1(a)=T$ on the domain of the canonical
projection $S^3\rightarrow\CPO$. By definition, the line bundle
$\gamma_1$ pulls back to
\[
S^3\times_T\bC\longrightarrow Q_1,
\]
where $T$ acts on $\bC$ by $w\mapsto u_1^{-1}w$.

For any $k\geq 1$ we assume that $\phi_k$ has been constructed with
the stated properties. So $\phi_k^*\gamma(a_k)$ is given by
\[
(S^3)^k\times_{T^k(a)}\bC\longrightarrow Q_k,
\]
where $T^k(a)$ acts on $\bC$ by $w\mapsto u_1^{-a(1,k+1)}\dots
u_k^{-a(k,k+1)}w$. It follows that the projectivisation ${\it
CP}(\phi_k^*(\bC\oplus\gamma(a_k)))$ coincides with
$Q_{k+1}$, and we define $\phi_{k+1}$ to be the resultant bundle map
to ${\it CP}(\bC\oplus\gamma(a_k))$. Then $\phi_{k+1}^*\gamma_j$ takes
the the required form for $1\leq j \leq k+1$.
\end{proof}

Form this point on we shall treat $Q_k$ and $M^k$ as
interchangeable, relating their properties by $\phi_k$ as necessary.
For example, the sections
\[
r_k,\;\widetilde{r}_k\colon Q_{k-1}\longrightarrow Q_k
\]
are induced by the inclusions of the respective subspaces
$(S^3)^{k-1}\times(1,0)$ and $(S^3)^{k-1}\times(0,1)$ of $(S^3)^k$,
using the notation of \eqref{sthreek}.

Following \cite{bura:tst}, we write the facets of $I^k$ as
$C_h^\varepsilon$, where $1\leq h\leq k$ and $\varepsilon$ is $0$ or
$1$. Thus $C_h^\varepsilon$ is the $(k-1)$--cube
$I^{h-1}\times(\varepsilon,1-\varepsilon)\times I^{k-h}$ in
$\bR_{\scriptscriptstyle\geqslant}^{2k}$. Every facet lifts to a
codimension--$2$ submanifold of $Q_k$, with normal $2$--plane bundle
$\nu_h^\varepsilon$. This is oriented if and only if the corresponding
isotropy subcircle $T(C_h^\varepsilon)<T^k_a$ is oriented, since
$T(C_h^\varepsilon)$ acts on the normal fibres. An {\it
omniorientation\/} of $Q_k$ is a choice of orientation for every
$\nu_h^\varepsilon$; there are therefore $2^{2k}$ omniorientations in
all, and each is preserved by the action of $T^k_a$.

The Pontryagin-Thom collapse maps $Q_k\rightarrow
T(\nu_h^\varepsilon)$ determine $2$--plane {\it facial bundles}
$\rho_h^\varepsilon$ over $Q_k$. Moreover, an orientation of
$\nu_h^\varepsilon$ determines, and is determined by, an orientation
of $\rho_h^\varepsilon$, for every $1\leq h\leq k$. An omniorientation
of $Q_k$ therefore identifies each of the $\rho_h^\varepsilon$ as
complex line bundles, and reversing any of the constituent
orientations induces complex conjugation on the corresponding line
bundle.
 
As explained in \cite{bura:tst}, there is a canonical isomorphism
\begin{equation}\label{taumkiso}
\tau(Q_k)\oplus\bR^{2k}\;\cong\;
\bigoplus_{h=1}^k\rho_h^0\oplus\rho_h^1
\end{equation}
of real $4k$--bundles. Every omniorientation therefore invests the
right-hand side with a complex structure, so that \eqref{taumkiso}
defines a corresponding stably complex structure on $Q_k$. These
structures play an interesting part in complex cobordism theory, and
we shall consider their enumeration in Section \ref{stcost}. As we
shall see, they include \eqref{stanmk}.

In \cite{grka:btc}, Grossberg and Karshon use a noncompact version of
\eqref{torusta} to describe Bott towers as complex manifolds. Given a
list $c=(c_1,\dots,c_{k-1})$ of integral sequences, they construct
$N_k$ as the quotient of $(\bC^2\setminus 0)^k$ by a $k$-fold
algebraic torus $\bC_\times^k(c)$, under the action
\begin{equation}\label{gkdef}
\begin{split}
&(w_1,\dots,w_k)\cdot(y_1,z_1;\dots;y_k,z_k)\;=\\
(y_1w_1,z_1w_1;y_2w_2,&w_1^{c(1,2)}z_2w_2;\dots;
y_kw_k,w_1^{c(1,k)}w_2^{c(2,k)}\dots w_{k-1}^{c(k-1,k)}z_kw_k).
\end{split}
\end{equation}
As complex manifolds, $N_k$ coincides with $Q_k$, where the latter is
determined by the list $a=-c$ (for which $a(i,j)=-c(i,j)$ for all $1\leq
i<j\leq k$). The corresponding structure on $M^k$ is that of the
projective form, introduced in Section \ref{boto}. These observations
are used in \cite{ci:phd} to relate the quotient cube $I^k$ to the
smooth fan determining $Q_k$.

Note that Grossberg and Karshon's construction yields the bounded flag
manifolds $B_k$ when $c_j=(0,\dots,0,-1)$, where $1\leq j\leq k-1$.

By mimicing the standard analysis for $\CP^k$ \cite{mist:cc}, we
deduce that the corresponding complex tangent bundle $\tau_\bC(Q_k)$
admits a canonical isomorphism
\begin{equation}\label{cancplx}
\tau_\bC(Q_k)\oplus\bC^k\;\cong\;
(\bC^2\setminus 0)^k\times_{\bC_\times^k(a)}\bC^{2k},
\end{equation}
where $\bC_\times^k(a)$ acts on $\bC^{2k}$ by extending \eqref{gkdef}.
The right-hand side splits as the sum of $2k$ complex line bundles,
where $\bC_\times^k(a)$ acts on $\bC$ by
\[
y_h\longmapsto y_hw_h\spandsp z_h\longmapsto
w_1^{-a(1,h)}\dots w_{h-1}^{-a(h-1,h)}z_hw_h,
\]
for $1\leq h\leq k$. Proposition \ref{qisproj} identifies these
bundles as $\overline{\gamma}_h$ and
$\overline{\gamma}_h\otimes\gamma(a_{h-1})$ respectively. So we may
rewrite \eqref{cancplx} as
\begin{equation}\label{canomni}
\tau_\bC(Q_k)\oplus\bC^k\;\cong\;
\bigoplus_{h=1}^k\overline{\gamma}_h\oplus
\big(\overline{\gamma}_h\otimes\gamma(a_{h-1})\big).
\end{equation}

The derivation of \eqref{taumkiso} yields isomorphisms
$\rho^0_h\cong\gamma_h$ and
$\rho^1_h\cong\overline{\gamma}_h\otimes\gamma(a_{h-1})$ of real
$2$--plane bundles. It follows that the stably complex structure
\eqref{canomni} arises from an omniorietation of $Q_k$. The structures
induced by the remaining $2^{2k}-1$ omniorientations may then be
obtained by replacing appropropriate line bundles by their complex
conjugates on the right-hand side of \eqref{canomni}. We use this
procedure to establish \eqref{omnidiffs} below.

%
%
%
%
%
%
%
%
%

\section{$KO$-Theory of Stages 1 and 2}\label{rekstontw}

The $KO$-theory of toric manifolds is considerably more subtle than
its complex counterpart, and is rarely free over the
coefficients. Bahri and Bendersky \cite{babe:ktt} have obtained
interesting results using the Adams Spectral Sequence, although their
calculations are mainly additive and make little reference to the
geometry of vector bundles. Our goal is to describe $\KO^*(M^k)$ as a
$\KO_*$-algebra for several families of Bott towers, in terms of the
bundles that we have introduced above. We also wish to understand the
complexification homomorphism, for application to stably complex
structures and cobordism classes in Section \ref{stcost}. Here we
focus on $M^1$ and $M^2$, which act as base cases for inductive
calculation and are useful for establishing notation.

It is convenient to denote the coefficient ring by
\[
\KO_*=\bZ[e,x,y]/(2e,e^3,ex,4x^2-y),
\]
where $e$, $x$, and $y$ are represented by the real Hopf line bundle
over $S^1$, the symplectic Hopf line bundle over $S^4$, and the
canonical bundle over $S^8$ respectively \cite{ka:kti}. We recall that
$\KO^*(S^n)$ is a free $\KO_*$-module on the single generator
$s_n^{KO}=s_n\in\KO^n(S^n)$, such that $s_n^2=0$ for each $n\geq 0$.

We appeal repeatedly to Bott's exact sequence
\begin{equation}\label{realexse}
\dots\longrightarrow\KO^{*-1}(X)\stackrel{\cdot e}{\longrightarrow}
\KO^{*-2}(X)\stackrel{\chi}{\longrightarrow}K^*(X)
\stackrel{r}{\longrightarrow}\KO^*(X)\longrightarrow\dots,
\end{equation}
which links real and complex $K$-theory through the realification
homomorphism $r$. Here, $\cdot e$ denotes multiplication by $e$, and
$\chi$ is defined by composing complexification $c$ with
multiplication by $z^{-1}$. For any element $g$ of $K^*(X)$, the
difference $g-\overline{g}$ lies in the kernel of $r$, and hence in
the image of $\chi$. Moreover,
\begin{equation}\label{cr}
c(r(g))=g+\overline{g},\sands\chi(r(zg))=g-\overline{g}.
\end{equation}
On the other hand, $r(c(h)=2h$ for any $h$ in $\KO^*(X)$. It is
important to remember that $c$ is multiplicative, whereas $r$ is not.

As in Fujii \cite{fu:kgp}, we define elements $u_i$ in
$\KO^{-2i}(\CPN)$ by $u_i=r(z^{i+1}u(n))$ for any integer $i$, where
$u(n)\in K^2(\CPN)$ arises in \eqref{ecpn}; as a ring, $\KO^*(\CPN)$
may then be described in terms of the $u_i$. When $n=2$, Fujii's
computations stretch to an isomorphism
\begin{equation}\label{kocpt}
\KO^*(\CPTp)\;\cong\;\KO_*[u_i:i\in\bZ]/F^2
\end{equation}
of $\KO_*$-algebras, where $F^2$ is the ideal
\[
\big(eu_i,\;xu_i-2u_{i+2},\;u_iu_{2j},\;
u_{2i+1}u_{2j-1}-4u_{2(i+j)}: \text{all $i$, $j$}\big).
\]
The relations show that $\KO^*(\CPT)$ is free of additive torsion, and
that $yu_i=u_{i+4}$ for all $i$; it therefore suffices to use $u_0$,
$u_1$, $u_2$, and $u_3$, as in \cite{fu:kgp}, but we retain the other
$u_i$ for notational convenience. We note that \eqref{kocpt} actually
defines a free $K_*$-module on a single generator $u_i$, where $zu_i$
is given by $u_{i+1}$ for any $i$.  This is equivalent to Wood's
well-known result \cite{wo:ktc} that $\KO\wedge\CPT$ is homotopy
equivalent to $K$.

Further computations lead to an isomorphism
\begin{equation}\label{kocpi}
\KO^*(\CPIp)\;\cong\;\KO_*[[u_i:i\in\bZ]]/F^\infty
\end{equation}
of $\KO_*$-algebras, where $F^\infty$ is the ideal
\[
\big(eu_i,\;xu_i-2u_{i+2},\;u_iu_j-u_{i-2}u_{j+2},\;
u_{2i+1}u_{2j-1}-(u_0+4)u_{2(i+j)}: \text{all $i$, $j$}\big).
\]
So $\KO^{2n}(\CPI)$ is torsion-free, and isomorphic to
$u_{-n}\bZ[[u_0]]$ for any integer $n$, whereas $\KO^{2n+1}(\CPI)$
is zero. For any complex line bundle $\gamma$ over a $2$--generated
complex $X$, it is convenient to interpret the pull-back of $u_i$
along the classifying map of $\gamma$ as a characteristic class
$u_i(\gamma)$ in $\KO^{-2i}(X)$.

It follows from \eqref{kocpi} that $\KO_*(\CPI)$ is torsion free,
and that $\KO\wedge\CPI$ is homotopy equivalent to the wedge
$\KO\wedge\big(\bigvee_{k\geq 0}\varSigma^{4k}\CPT\big)$. This
equivalence may also be deduced from the fact that a vector bundle
is $\KO$-orientable precisely when it is {\it Spin\/}
\cite{atbosh:cm}.

We consider $\zeta^2$ over $\CPI$, which is universal for complex
line bundles with {\it Spin\/}-structure, and utilise the Thom class
$t^K$ of Lemma \ref{twogen} in $K^2(T(\zeta^2))$.
\begin{lem}\label{tzetasq}
There is a unique element $\tquare$ in $\KO^2(T(\zeta^2))$ whose
complexification is given by $c(\tquare)=\overline{\zeta}t^K$; it is
a Thom class, and satisfies $t_{\quare}^2=u_{-1}\tquare$ in
$\KO^4(T(\zeta^2))$.
\end{lem}
\begin{proof}
The existence of a Thom isomorphism
$\KO^{*-2}(\CPIp)\cong\KO^*(T(\zeta^2))$ confirms that
$\KO^{2n}(T(\zeta^2))$ is torsion free for $n\not\equiv 3\bmod 4$. So
\eqref{realexse} reduces to a short exact sequence
\[
0\longrightarrow\KO^{2n}(T(\zeta^2))\stackrel{\chi}{\longrightarrow}
K^{2n+2}(T(\zeta^2))\stackrel{r}{\longrightarrow}\KO^{2n+2}(T(\zeta^2))
\longrightarrow 0,
\]
for $n=1$ and $2$; thus $c$ is monic, and if $\tquare$ exists, it is
unique.

The construction of $t^K$ implies that
$\overline{t^K}=\overline{\zeta}^2t^K$, so that
\[
c\cdot r(z^{-1}\overline{\zeta}t^K)=
z^{-1}(\overline{\zeta}t^K-\zeta\overline{t^K})=0;
\]
hence $r(z^{-1}\overline{\zeta}t^K)=0$, and $\tquare$ exists as
required. It is a Thom class because $\overline{\zeta}t^K$ is a Thom
class and $c$ is a map of ring spectra. Moreover,
$(t^K)^2=z^{-1}(\zeta^2-1)t^K$, whence
\[
c(t_{\quare}^2)=z^{-1}(1-\overline{\zeta}^2)t^K=
z^{-1}(\zeta-\overline{\zeta})c(\tquare)=c(u_{-1}\tquare).
\]
Thus $t_{\quare}^2=u_{-1}\tquare$ in $\KO^4(T(\zeta^2))$.
\end{proof}

The calculation of $\KO^*(T(a))$ depends on the parity of $a$. When
$a=2b$ is even, $\zeta(1)^a$ is {\it Spin}$(2)$-bundle, and is the
pull-back of the universal example along the map $\CPO\rightarrow\CPI$
of degree $b$; thus $\tquare$ pulls back to a Thom class $t$ in
$\KO^2(T(a))$.
\begin{prop}\label{kotl}
When $a$ is even, $\KO^*(T(a))$ is isomorphic to
\[
\KO_*[s_2,t]\big/(s_2^2,t^2-as_2t)
\]
as $\KO_*$-algebras. When $a$ is odd, there are elements $m_i$ in
$\KO^{-2i}(T(a))$ such that $\KO^*(T(a))$ is isomorphic to
\[
\KO_*[m_i:i\in\bZ]/F(a,m)
\]
as $\KO_*$-algebras, where $F(a,m)$ is the ideal
\[
\big(em_i,\;xm_i-2m_{i+2},\;m_im_{2j},\;
m_{2i+1}m_{2j-1}-4am_{2(i+j)}:
\text{all $i$, $j$}\big).
\]
\end{prop}
\begin{proof}
When $a$ is even, the Thom isomorphism identifies $\KO^*(T(a))$ with
the free $\KO_*$-module on generators $t$ and $s_2t$. It therefore
remains to evaluate $t^2$ in $\KO^4(T(a))$. But
$u_{-1}(\zeta(1)^b)=r(z^{-1}(\zeta(1)^b-1))$ in $\KO^2(\CPO)$, so
$t^2=br(z^{-1}(\zeta(1)-1))t=as_2t$, as required.

When $a=2b+1$ is odd, $\zeta(1)^a$ is no longer $\KO$-orientable. We
proceed by comparing the $\KO$-theory of the cofibre sequences of
$S^2\cup_{a\eta}e^4$ and $\CPT$, using the map $f(a)\colon
T(a)\rightarrow\CPT$ which classifies $\zeta(1)^a$. We define $m_i$ as
$r(z^{i+1}(1-bzs_2^K)t^K)$ when $i$ is even, and $r(z^{i+1}t^K)$ when
$i$ is odd. The action of $f(a)^*$ then yields the algebra structure,
by appeal to \eqref{kocpt}; alternatively, we may apply
complexification.
\end{proof}

A few observations are in order. Firstly, when $a$ is even the
suspension of $a\eta$ is null homotopic, so that $\varSigma T(a)$ is
homotopy equivalent to $S^3\vee S^5$; equivalently, the \mbox{\it
SO}$(3)$-bundle $\bR\oplus\zeta(1)^a$ is trivial. Secondly, the
relations of Proposition \ref{kotl} imply that $t^3=0$. Thirdly, the
action of $f(a)^*$ is computed from \eqref{realexse}, and is given by
\begin{equation}\label{faeff}
f(a)^*(u_i)=
\begin{cases}
(2+be^2s_2)t&i=-1\\
bxs_2t&i=0\\
xt&i=1\\
ays_2t&i=2
\end{cases}
\;\;\text{and}\;\;
f(a)^*(u_i)=
\begin{cases}
am_i&i\equiv 0(2)\\
m_i&i\equiv1(2)
\end{cases}\;
\end{equation}
for $a=2b$ and $2b+1$ respectively. Fourthly, when $a$ is odd, the
generators $m_i$ may be defined more systematically as
$r(z^{i+1}\overline{\zeta}(1)^bt^K)$; this description is central to
Theorem \ref{kobtto} below.

Proposition \ref{kotl} shows that $\KO^*(T(a))$ is free over $\KO_*$
when $a$ is even, and over $K_*$ when $a$ is odd. It may be
interpreted in terms of spectra as providing homotopy equivalences
\begin{equation}\label{kothom}
\KO\wedge T(2b)\;\simeq\;\KO\wedge
\left(S^2\vee S^4\right)\sands
\KO\wedge T(2b+1)\;\simeq\;\KO\wedge\CPT.
\end{equation}

We may now proceed to $M^2$ via Proposition \ref{btsplit}, which
ensures that there is an additive isomorphism
\begin{equation}\label{komtadd}
\KO^*(M^2)\;\cong\;\KO^*(S^2)\oplus\KO^*(T(a))
\end{equation}
of $\KO_*$-modules. It remains to describe the products in
$\KO^*(M^2)$. To prepare for our eventual notation, we write
$p_2^*s_2$ as $d_1$ in $\KO^2(M^2)$ and $q_2^*t$ as $d_2$ in
$\KO^2(M^2)$, when $a$ is even; when $a$ is odd, we write $q_2^*m_i$
as $n_i$ in $\KO^{-2i}(M^2)$, for all $i$.
\begin{prop}\label{komt}
When $a$ is even, $\KO^*(M^2_+)$ is isomorphic to
\[
\KO_*[d_1,d_2]\big/(d_1^2,\,d_2^2-ad_1d_2)
\]
as $\KO_*$-algebras; when $a$ is odd, it is isomorphic to
\[
\KO_*[d_1,n_i:i\in\bZ]\,\big/\big(F(a,n),\,d_1^2,\,d_1n_{2i},
d_1n_{2i+1}-2n_{2i}\big).
\]
\end{prop}
\begin{proof}
It suffices to combine Proposition \ref{kotl} with \eqref{komtadd}.
When $a$ is odd, the extra relations follow by applying
complexification, and noting that $n_i$ restricts to $0$ on $M^1$ for
all $i$.
\end{proof}

The following corollary is immediate, and helps us to enumerate stably
complex structures on $M^2$ in Section \ref{stcost}.
\begin{cor}\label{remt}
In both cases, $\KO^{-2}(M^2)$ is isomorphic to $\bZ^2$ as abelian
groups; bases are given by $\{xd_1, xd_2\}$ when $a$ is even, and
$\{xd_1,n_1\}$ when $a$ is odd.
\end{cor}

The isomorphisms of \eqref{mtisos} extend to $\KO^*(M^2)$, and may
be described in terms of \eqref{kothom} and Proposition \ref{komt}.

%
%
%
%
%
%
%
%
%

\section{$KO$-Theory of Bott Towers}\label{rekboto}

We now return to the Bott tower $(M^k:k\leq n)$, determined by the
list $a=(a_1,\dots,a_{n-1})$, and study inductive procedures for
computing the $\KO_*$-algebra structure of $\KO^*(M^k)$ in favourable
cases.

The work of Bahri and Bendersky \cite{babe:ktt} identifies the effect
of smashing $M^k$ with the spectrum $\KO$, and leads to a homotopy
equivalence
\begin{equation}\label{babehe}
\KO\wedge N^{2n}_+\;\simeq\;\KO\wedge
\bigvee_{p,q=0}^{n,n-2}\left(\bigvee^{\alpha_p}S^{2p}
\vee\bigvee^{\beta_q}\varSigma^{2q}\CPT \right)
\end{equation}
for any toric manifold $N^{2n}$. The {\it BB-numbers} $\alpha_p$ and
$\beta_q$ enumerate the summands for each $p$ and $q$
respectively. \baabe\ prove that their numbers are determined by the
structure of $H^*(N^{2n};\bF_2)$ over $\mathcal{A}(1)$, the subalgebra
of the Steenrod algebra generated by $\Sq{1}$ and $\Sq{2}$. Two types
of $\mathcal{A}(1)$-module are involved; the first is
$\varSigma^{2p}\mathcal{M}_1$, with one $2p$--dimensional generator on
which $\Sq{1}$ and $\Sq{2}$ act trivially, and the second is
$\varSigma^{2q}\mathcal{M}_2$, with one $2q$--dimensional generator
$x$ such that $\Sq{1}x=0$ and $\Sq{2}x\neq 0$. Then
$H^*(N^{2n};\bF_2)$ decomposes as a direct sum of these two types; the
number of summands $\varSigma^{2p}\mathcal{M}_1$ is $\alpha_p$, and
the number of summands $\varSigma^{2q}\mathcal{M}_2$ is $\beta_q$.

The additive part of our calculations recover \eqref{babehe} for two
particular families of Bott towers, and provide representative bundles
for the generators of $\KO^*(M^k)$ as a geometrical bonus. We also
point out how the BB-numbers depend on the parity of the entries in
$a$. Our families actually illustrate the extreme cases, which range
from $\beta_q=0$ for all $q$, to $\alpha_p=0$ for all $p>1$.

We begin by reverting to the notation of Section \ref{twgeco}, and
consider the complex line bundle
$\gamma=\gamma_1^{a(1)}\otimes\dots\otimes\gamma_m^{a(m)}$ over the
2--generated complex $X$.

When $a(j)=2b(j)$ is even for all $1\leq j\leq m$, we write
$\gamma_1^{b(1)}\otimes\dots\otimes\gamma_m^{b(m)}$ as
$\gamma^{1/2}$. So $\gamma$ is {\it Spin}$(2)$, and is obtained by
pulling the universal example of Lemma \ref{tzetasq} back along the
classifying map for $\gamma^{1/2}$. In particular, we obtain a Thom
class $t\in\KO^2(T(\gamma))$; it satisfies
$t^2=u_{-1}(\gamma^{1/2})t$, where $u_{-1}(\gamma^{1/2})=
r(z^{-1}(\gamma^{1/2}-1))$ in $KO^2(X)$, and
\[
c(t)=\prod_{j\leq m}(1+zg_j)^{-b(j)}t^K
\]
in $K^0(T(\gamma))$.
\begin{prop}\label{tsq}
The $\KO_*$-algebra $\KO^*(Y_+)$ is a free module over $\KO^*(X_+)$ on
generators $1$ and $d_{m+1}$, which have dimensions $0$ and $2$
respectively; the multiplicative structure is determined by the single
relation
\begin{equation}\label{dsqrel}
d_{m+1}^2=r\big(z^{-1}\big(\prod_{j\leq
m}(1+zg_j)^{b(j)}-1\big)\big)d_{m+1},
\end{equation}
and $d_{m+1}$ restricts to a generator on the fibre $S^2\subset Y$.
\end{prop}
\begin{proof}
We repeat the arguments of Lemma \ref{twogen}(2) with $q^*t=d_{m+1}$
in $\KO^2(Y)$, and apply the remarks above.
\end{proof}

It is sometimes preferable to leave \eqref{dsqrel} in the form
$d_{m+1}^2=u_{-1}(\gamma^{1/2})d_{m+1}$, and aim to express
$u_{-1}(\gamma^{1/2})$ as a polynomial in the elements $r(z^ig_j)$.
This does not follow automatically from \eqref{dsqrel}, because $r$ is
not multiplicative. The simplest example is $X=S^2$, where
$\gamma^{1/2}$ is given by $\zeta^{b(1)}$ and $u_{-1}(\zeta^{b(1)})$
reduces to $2b(1)s_2$ in $\KO^2(S^2)$. We then recover the first part
of Proposition \ref{komt}.

If one or more of the integers $a(j)$ is odd, the situation is less
amenable. For our current purposes, it is enough to recall that
$T(\gamma)$ admits a canonical complex line bundle $\lambda$ over
$T(\gamma)$, defined by $v^H(\lambda)=t^H$. So $t^K$ is represented by
$z^{-1}(\lambda-1)$ in $K^2(T\gamma)$. The classes $u_i(\lambda)$ in
$\KO^{-2i}(T(\gamma))$ play a major r\^ole in describing $\KO^*(Y_+)$.

Our main structure theorems refer to two particular families of Bott
towers. They are the {\it totally even\/} towers, for which the
integers $a(i,j)=2b(i,j)$ are even for all values of $1\leq i<j\leq
n$, and the {\it terminally odd\/} towers, for which the integers
$a(j-1,j)=2c(j)+1$ are odd for every $1\leq j\leq n$. It is possible
to deal with other cases by combining the two approaches.

\begin{thm}\label{kobtte}
For any totally even Bott tower $(M^k:k\leq n)$, the $\KO_*$-algebra
$\KO^*(M^k_+)$ is isomorphic to $\KO_*[d_1,\dots,d_k]/J^{te}_k$,
where $J^{te}_k$ denotes the ideal
\[
\Big(d_j^2-r\Big(z^{-1}
\Big(\prod_{i<j}(1+zg_i)^{b(i,j)}-1\Big)\Big)d_j: 1\leq j\leq k\Big);
\]
for each $1\leq k\leq n$, the homotopy equivalence $h_k$ induces the
$\KO_*$-module isomorphism
\[
\KO^*(M^k)\cong\langle\dlo\rangle\oplus\dots\oplus
\langle\dlk\rangle,
\]
where $\langle\dlj\rangle$ denotes the free submodule generated by
those monomials $d_R$ for which $R\subseteq[j]$ and $j\in R$.
\end{thm}
\begin{proof}
In this case the proof of Theorem \ref{cohomk} adapts directly, since
all the relevant $\KO_*$-modules are free.
\end{proof}

As before, it may be preferable to rewrite the relations of $J^{te}_k$
as
\begin{equation}\label{alttoev}
d_j^2\;=\;u_{-1}(\gamma_1^{b(1,j)}\otimes\dots
\otimes\gamma_{j-1}^{b(j-1,j)})\,d_j,
\end{equation}
and calculate $u_{-1}(\gamma_1^{b(1,j)}\otimes\dots
\otimes\gamma_{j-1}^{b(j-1,j)})$ as a polynomial in $d_1$, \dots,
$d_{j-1}$ for each $1\leq j\leq k$. Amongst other formulae in
$\KO^*(M^k)$, this approach yields
\[
d_j^{j+1}=0\sands
d_j^2=\big(a(1,j)d_1+\dots+a(j-1,j)d_{j-1}\big)d_j\;\;\text{modulo
$P_*$},
\]
where $P_*$ denotes the ideal generated by triple products.

Calculations for terminally odd towers are more intricate, and we
begin with the additive structure. It is convenient to index the
generators by finite sets $R$ of positive integers. For every such $R$,
we construct $R^+$ by adding $1$ to each element, and $1;R^+$ by
adjoining the integer $1$ to the result. We obtain the coproduct
decomposition
\begin{equation}\label{coppow}
2^{[k-1]}\stackrel{e_1}{\llongrightarrow}2^{[k]}
\stackrel{e_2}{\llongleftarrow}2^{[k-1]},
\end{equation}
of power sets, where $e_1(R)=R^+$ and $e_2(R)=1;R^+$. Given
$R\subseteq [k-2]$ for $k\geq 2$, we construct $R;k\subset [k]$ by
adjoining the integer $k$.

So far as complex $K$-theory is concerned, we may apply this notation
to the ladder \eqref{ladder}. The elements $g_Rt^K_k$ in
$K^*(T(a_{k-1}))$ are of two types; those for which $R$ takes the form
$S^+$ for some $S\subseteq [k-3]$, so that
$i^*(g_Rt^K_k)=g'_S(t')^K_{k-1}$, and those for which $R$ takes the
form $1;S^+$, so that $f^*(\varSigma^2g'_S(t')^K_{k-1})=g_Rt^K_k$.
The decomposition \eqref{coppow} then corresponds to the splitting
\eqref{split1}. Of course, $q_k^*(g_Rt^K_k)=g_{R;k}$ in $K^*(M^k)$.

We may now construct the elements we need in $\KO$-theory. For every
integer $i$, we define
\[
m(R;k)_i\;=\;r\big(z^{i+1}\hpm
\overline{\gamma}_{k-1}^{\,b(k)}\hpm g_Rt^K_k\big)
\]
in $\KO^{2(|R|-i)}(T(a_{k-1}))$, as $R$ ranges over subsets of
$[k-2]$, and
\[
n(R;j)_i\;=\;r\big(z^{i+1}\hpm
\overline{\gamma}_{j-1}^{\,b(j)}\hpm g_{R;j}\big)
\]
in $\KO^{2(|R|-i)}(M^k)$, as $R$ ranges over subsets of $[j-2]$,
with $2\leq j\leq k$. Thus $q_k^*m(R;k)_i=n(R;k)_i$ for every
$R\subseteq [k-2]$.

\begin{thm}\label{kobtto}
For any terminally odd Bott tower $(M^k:k\leq n)$, the
$\KO_*$-module $\KO^*(M^k_+)$ is generated by the elements
\[
\big\{d_1,\;n(R;j)_i:2\leq j\leq k\big\},
\]
where $R$ ranges over the subsets of $[j-2]$ and $i\in\bZ$; the
submodule of relations is generated by 
\[
\big\{en(R;j)_i,\; xn(R;j)_i-2n(R;j)_{i+2}\big\}
\]
for all $R$, $j$ and $i$.
\end{thm}
\begin{proof}
We proceed inductively, using the commutative ladder \eqref{ladder}.
We assume that the result holds for terminally odd towers of height
$\leq n-1$, where $n\geq 2$, and consider $(M^k:k\leq n)$,
determined by a list $(a_1,\dots,a_{n-1})$. The tower $((M')^k:k\leq
n-1)$ is determined by the list $(a_1',\dots,a_{n-2}')$, where
$a_{j-1}'$ is obtained from $a_j$ by deleting the first element; so
it is also terminally odd, and the inductive hypothesis applies.

We may therefore assume that $\KO^*(T(a_{k-2}'))$ is a free abelian
group, generated by the elements $m(S;k-1)'_i$ for $S\subseteq [k-3]$
and $0\leq i\leq 3$. So $\KO^*(\varSigma^2T(a_{k-2}'))$ is generated
by their double suspensions, and both groups are zero in odd
dimensions. Since $i^*\gamma_j=\gamma'_{j-1}$ in $\KO^0((M')^{k-2}$
for every $2\leq j\leq k-1$, it follows from Corollary \ref{decomptwo}
that $i^*m(S^+;k)_i=m(S;k-1)'_i$ in $\KO^0((M')^{k-2}$, and
$f^*(\varSigma^2m(S;k-1)'_i=m(1;S^+;k)_i$ in $\KO^*(T(a_k))$, for
every $S\subseteq [k-3]$ . Applying $\KO^*(-)$ to the ladder
yields
\[\begin{CD}
@<\delta\!\npm<<\KO^*(T(a_{k-2}'))\!
@<i^*<<\!\KO^*(T(a_{k-1}))\!@<f^*<<
\!\KO^*(\varSigma^2T(a_{k-2}'))\!\npm@<\delta<<
\\
@.@V(q_{k-1}')^*VV@Vq^*_kVV@VV\varSigma^2(q_{k-1}')^*V\\
@<<\delta<\!\npm\KO^*((M')^{k-1})\!@<<i^*<\KO^*(M^k)\!@<<f^*<
\!\KO^*(\varSigma^2(M')^{k-1}_+)\!\npm@<<\delta<
\end{CD}\!,
\]
ensuring that the upper coboundary maps $\delta$ are zero for $k\geq
2$, and that the upper sequence splits as abelian groups. So
$\KO^*(T(a_{k-1}))$ is also zero in odd dimensions, and generated by
the $m(S^+;k)_i$ and $m(1;S^+;k)_i$ in even dimensions; but these are
precisely the elements $m(R;k)_i$ for $R\subseteq [k-2]$. It follows
from Proposition \ref{btsplit} that $q_k^*$ injects
$\KO^*(T(a_{k-1}))$ into $\KO^*(M^k)$ as the summand generated by the
elements $n(R;k)_i$, for $R\subseteq [k-2]$. The abelian group
structure of $\KO^*(T(a_{k-1}))$ ensures that complexification is
monic, and therefore that
\[
en(R;k)_i=0\sands xn(R;k)_i=2n(R;k)_{i+2}
\]
for all $2\leq j\leq k$. The remainder of the additive structure then
follows from the inductive hypothesis. The base case $k=2$ is resolved
by Proposition \ref{komt}, with $m(\varnothing;2)_i=m_i$ and
$n(\varnothing;2)_i=n_i$ for all $i$.
\end{proof}

It follows from Theorem \ref{kobtto} that $\KO^*(M^k)$ is torsion
free, except for a single copy of $\bZ/2$ in each of the dimensions
$0$ and $1\bmod 8$, generated by $e^2y^id_1$ and $ey^id_1$
respectively. This generalises the results obtained for $B_k$ in
\cite{ci:phd}, and outlined in Example \ref{teodex} below. Theorem
\ref{kobtto} also implies that $yn(R;j)_i=n(R;j)_{i+4}$ in
$\KO^*(M^k)$, for any $R$, $j$, and $i$. We may therefore restrict the
choice of generators to $i=0$, $1$ $2$, and $3$, for example;
nevertheless, we usually allow $i$ to be arbitrary for notational
convenience.

In order to understand the multiplicative structure of $\KO^*(M^k)$, 
we need to evaluate products of the generators described in Theorem
\ref{kobtto}. 
\begin{prop}\label{bttoprods}
For any $R\subseteq [j-2]$ and $j>2$, we have that 
\[
d_1n(R;j)_i=
\begin{cases}
0&\text{if $1\in R$}\\
n(1;R;j)_i&\text{otherwise}
\end{cases}\;;
\]
for any $R'\subseteq [j'-2]$, we have that 
\begin{equation}\label{koprodform}
\begin{split}
n(R;j)_i\hpm n&(R';j')_{i'}\;=\;\\
&r\Big(z^{i+j+2}\hpm\overline{\gamma}_{j-1}^{\,b(j)}\hpm g_{R;j}
\Big(\overline{\gamma}_{j'-1}^{\,b(j')}g_{R';j'}+
(-1)^{j+1}\gamma_{j'-1}^{b(j')}\hpm\overline{g}_{R';j'}\Big)\Big).
\end{split}
\end{equation}
In particular, $n(R;j)_i\hpm n(R';j)_{i'}=0$ whenever $1\in R\cap R'$.
\end{prop}
\begin{proof}
Theorem \ref{kobtto} implies that complexification is monic, modulo
the summand $\KO^*(M^1)$. Since $n(R;j)_i$ restricts to $0$ in
$\KO^*(M^1)$ for every $R$, $j$, and $i$, it suffices to prove the
relations by applying $c$.

Now $c(d_i)=g_1$, and $c(n(R;j)_i)=
z^{i+1}(\overline{\gamma}^{\,b(j)}_{j-1}\hpm g_{R;j}+
(-1)^{i+1}\gamma^{b(j)}_{j-1}\hpm\overline{g}_{R;j})$ in
$K^*(M^k)$. Moreover, $g_1^2=0$, so $\overline{g}_1=g_1$ and the first
set of relations follows. The second set is proven similarly, by
noting that
\[
c\big(r(x)r(y)\big)=cr(x(y+\overline{y}))
\]
for any elements $x$ and $y$ in $\KO^*(M^k)$.  
\end{proof}
We would like to write \eqref{koprodform} as an explicit
$\KO_*$-linear combination of the generators $d_1$ and $n(R;j)_i$. In
principle, this may be achieved by using the expressions for
$\overline{g}_m$ and $g_m^2$ of \eqref{conjongj} and \eqref{eehak}
respectively; in practice, the calculations increase rapidly in
complexity. Examples \ref{toevex} and \ref{teodex} give a more
detailed glimpse of the difficulties which characterise the
multiplicative structures described in Theorem \ref{kobtte} and
Proposition \ref{bttoprods}. Related calculations will be presented in
\cite{do:phd}.

The following observations flow directly from Theorems
\ref{kobtte} and \ref{kobtto}.
\begin{cor}\label{kocor}
In the totally even case, the equivalence \eqref{babehe} reduces to
\[
\KO\wedge M^k_+\;\simeq\;\KO\wedge
\bigvee_{R\,\subseteq\,[k]}S^{2|R|}\,;
\]
thus $\alpha_p=\binom{k}{p}$ for all $1\leq p\leq k$, and $\beta_q=0$
for all $q$. In the terminally odd case, we have
\[
\KO\wedge M^k_+\;\simeq\;\KO\wedge
\left(S^2_+\vee\bigvee^{k-2}_{h=0\vphantom{|}}\;
\bigvee_{R\subseteq[h]}\varSigma^{2|R|}\CPT\right);
\]
thus $\alpha_p=0$ for $2\leq p\leq k$, and
$\beta_q=\sum_{h=q}^{k-2}\binom{h}{q}$ for all $0\leq q\leq k-2$.
\end{cor}
\begin{proof}
In the totally even case, Theorem \ref{kobtte} confirms that
$\KO^*(M^k_+)$ is additively generated over $\KO_*$ by the monomials
$d_R=\prod_Rg_j$, as $R$ ranges over the subsets of $[k]$.

In the terminally odd case, the torsion subgroup of $\KO^*(M^k)$
corresponds to the summand $\KO\wedge S^2$. The proof of Theorem
\ref{kobtto} combines with \eqref{kothom} to show that
\[
\KO\wedge T(a_{j-1})\;\simeq\;
\KO\wedge\left((S^2_+)^{\wedge(j-2)}\wedge\CPT\right)
\]
for all $1\leq j\leq k$, where the elements $n(R;j)_i$ correspond to
the summand $\varSigma^{2|R|}\CPT$ for every $R\subseteq [j-2]$. The
result now follows from Proposition \ref{btsplit}.
\end{proof}

Corollary \ref{kocor} illustrates the relationship between the
BB-numbers and entries in the list $a$. In the totally even case,
Proposition \ref{cohomk} confirms that every square is zero in
$H^*(M^k;\bF_2)$, so $\Sq{2}=0$; thus $\varSigma^{2q}\mathcal{M}_2$
cannot occur in its decomposition, and $\beta_q=0$ for all $1\leq
j\leq k$, as required. In the terminally odd case, we write the mod
$2$ reduction of the class $x_i$ as $x'_i$. Then Proposition
\ref{cohomk} confirms that $\Sq{2}x'_1=0$, and $\Sq{2}x'_j\equiv
x'_{j-1}x'_j$ modulo terms of the form $x'_ix'_j$ with $i\leq j-2$,
for every $2\leq j\leq k$. Thus $\alpha_1=1$. A simple inductive
calculation reveals that $H^{2q+2}(M^k;\bF_2)$ decomposes as
\[
\Sq{2}H^{2q}(M^k;\bF_2)\oplus H_{2q+2},
\]
where $H_{2q+2}$ is generated by all monomials of the form $x'_Rx'_j$
such that $R\subseteq [j-2]$ and $|R|=q$. Since $\Sq{2}$ is injective
on $H_{2q+2}$, it follows that $\alpha_p=0$ for $2\leq p\leq k$, and
$\beta_q=\sum_{h=q}^{k-2}\binom{h}{q}$ for all $q$, as required.

In order to illustrate these results, we discuss two examples.
\begin{exa}\label{toevex}
Let $(A_k:k\geq 0)$ be the totally even tower determined by the
integers $a(i,j)=0$ for $i\leq j-2$, and $a(j-1,j)=2$, for any $j\geq
1$. The relation \eqref{alttoev} reduces to
$d_j^2=u_{-1}(\gamma_{j-1})d_j$, so we have to compute
$u_{-1}(\gamma_{j-1})$ in $\KO^2(M^j)$; this follows inductively from
an understanding of the homomorphism
$f^*\colon\KO^*(\CPI)\rightarrow\KO^*(T(\zeta^2))$, where $f$ is the
map of Thom complexes classifying $\zeta^2$. To calculate $f^*$, we
extend the formulae of \eqref{faeff} in case $b=1$, and find
\begin{equation}\label{feff}
f^*(u_i)=
\begin{cases}
(2+u_0)\tquare &i=-1\\
(e^2+u_1)\tquare &i=0\\
(x+u_2)\tquare &i=1\\
u_3\tquare &i=2
\end{cases}.
\end{equation}
We deduce that $u_{-1}(\gamma_{j-1})$ is given by 
\[
\sum_{s=0}^{\lfloor\frac{j-1}{4}\rfloor}2y^sd_{j-1}\dots d_{j-4s}+
\sum_{s=0}^{\lfloor\frac{j-2}{4}\rfloor}e^2y^sd_{j-1}\dots d_{j-4s-1}+
\sum_{s=0}^{\lfloor\frac{j-4}{4}\rfloor}xy^sd_{j-1}\dots d_{j-4s-3}.
\]
\end{exa}
\begin{exa}\label{teodex}
Let $(B_k:k\geq 0)$ denote the terminally odd tower of bounded flag
manifolds, determined by integers $a(i,j)=0$ for $i\leq j-2$ and
$a(j-1,j)=1$, for all $j\geq 1$. Then each $b(j)$ is zero, and the
generators $n(R;j)_i$ are defined by $r(z^{i+1}g_{R;j})$ for every
$R\subseteq [j-2]$. Products of the form $n(R;j)_i\cdot n(R';j')_{i'}$ 
are given by 
\[
r\big(z^{i+j+2}(g_{R;j}\hpm g_{R';j'}+
(-1)^{j+1}g_{R;j}\hpm\overline{g}_{R';j'})\big),
\]
and are evaluated using the formulae
\[
g_m^2=
\Bigg(\sum_{\varnothing\neq S\subseteq [m-1]}z^{|S|-1}g_S\Bigg)g_m
\spandsp\overline{g}_m=g_m/(1+zg_m)
\] 
in $K^*(M^k)$, for every $1\leq m\leq k$.
\end{exa}

We may combine Theorems \ref{kobtte} and \ref{kobtto} to identify
$\KO^{-2}(M^k)$. As explained in Section \ref{stcost}, these groups
classify the stably almost complex structures on $M^k$.  
\begin{thm}\label{komintwo}
If the tower is totally even, then $\KO^{-2}(M^k)$ is isomorphic to
\[
\Big(\bigoplus_{|R|\equiv 1,-1(4)}\bZ\Big)\;\oplus\;
\Big(\bigoplus_{|R|\equiv 0(4)}\bZ/2\Big),
\]
where $R\subseteq[k]$; a basis is given by 
\[
\big\{xy^{(|R|-1)/4}d_R,\,y^{(|R|+1)/4}d_R,\,e^2y^{|R|/4}d_R\big\}.
\]
If the tower is terminally odd, then $\KO^{-2}(M^k)$ is isomorphic to
$\bZ^{2^{k-1}}$; a basis is given by $\big\{xd_1,n(R;j)_i\big\}$,
where $R\subseteq [j-2]$ for $2\leq j\leq k$ and $i=|R|+1$.
\end{thm}

%
%
%
%
%
%
%
%
%

\section{Stably Complex Structures}\label{stcost}

By way of conclusion, we apply our results to the study of stably
complex structures on certain families of Bott towers. We consider the
enumeration of those which arise from omniorientations, and discuss
two particular special cases; those which restrict to almost complex
structures, and those which are null-cobordant in $\osu$. We summarise
the appropriate definitions in order to establish our notation.

Full details of the results for $(B_k:k\geq 0)$ in Theorems
\ref{omnisscs}, \ref{omniacs} and \ref{omnibds} are provided in
\cite{ci:sac}. 

We write $\BU$ and $\BO$ respectively for the classifying spaces of
the infinite unitary and orthogonal groups, and let 
$r\colon\BU\rightarrow\BSO\subset\BO$ denote a specific choice of
realification. The resulting maps
\[
\SO/U\stackrel{f}{\longrightarrow}\BU\stackrel{r}{\longrightarrow}\BO
\]
induce the $K$-theory exact sequence \eqref{realexse} for connected
spaces $X$. Given a smooth oriented manifold $N$, we assume that the
stable tangent bundle is represented by a map $\tau^S\colon
N\rightarrow\BSO$, which we fix henceforth. A complex structure on
$\tau^S$ is given by a lift $\tau$ to {\it BU}, and is known as a
stably complex structure, or {\it U-structure}, on $N$; it therefore
consists of a factorisation $\tau^S=r\cdot\tau$. We deem two
$U$-structures $\tau$ and $\tau'$ to be equivalent, or {\it
homotopic}, whenever they are homotopic through lifts of $\tau^S$.
Once $\tau$ is chosen, it leads to a complementary lift of the stable
normal bundle $\nu^S$ of $N$, and conversely; this correspondence
preserves homotopy classes.

If we begin with the opposite orientation for $N$, we obtain a second
set of $U$-structures and homotopy classes. They are distinct from
those described above, but correspond to them bijectively. 
 
An {\it almost complex\/} structure on $N$ is given by a complex
structure on the tangent bundle $\tau(N)$, and determines a compatible
orientation. When $N$ is a complex manifold, it therefore admits a
corresponding almost complex structure, which stabilises to the {\it
underlying} $U$-structure $\tau_\bC$. An arbitrary $U$-structure need
not, of course, destabilise to $\tau(N)$, just as an almost complex
structure need not be integrable. Henceforth, we will deal only with
complex connected $N$, oriented compatibly, and will take $\tau_\bC$
to be the distinguished $U$-structure. As explained in
\cite{raswta:gsg}, we may then define a bijection between
$\KO^{-2}(N)$ and the homotopy classes of $U$-structures on $N$. To
each $\varDelta\in\KO^{-2}(N)$ there corresponds a homotopy class of
complex structures on the trivial bundle $\bR^{2L}$, for suitably
large $L$, and the bijection associates the $U$-structure
$\tau\letbe\tau_\bC\oplus\bR^{2L}$ to $\varDelta$. In other words,
$\varDelta(\tau,\tau_\bC)$ is the {\it difference element\/} of
$\tau$; its image under $\chi$ is represented by the virtual bundle
$\tau-\tau_\bC$ in $K^0(N)$.

So Theorem \ref{komintwo} identifies the totality of $U$-structures on
the Bott tower $(M^k:k\leq n)$. In the terminally odd case, $\chi$ is
monomorphic and the structures may be enumerated by identifying
$\tau_\bC$ as an element of $K^0(M^k)$, then varying $\tau-\tau_{\bC}$
over the image of $\chi$. This strategy was applied to the tower of
bounded flag manifolds $(B_k:0\leq k)$ in \cite{ci:phd}.

For more general purposes, it helps to follow the lead of Section
\ref{tost}, and define a complex structure on an arbitrary vector
bundle $\theta$ as an isomorphism $g$ from $\theta$ to a complex
vector bundle $\xi$. The action of $i$ on the fibres of $\theta$ is
given by conjugating its action on $\xi$ by $g$, and homotopy classes
of isomorphisms correspond to homotopy classes of complex
structures. An isomorphism of the form $\tau(N)\oplus\bR^m\cong\xi$
therefore specifies a $U$-structure on $N$; for example,
\eqref{canomni} defines the $U$-structure $\tau_\bC$ underlying the
projective form of $M^k$.

A second isomorphism $g'\colon\theta\cong\xi$ defines a second complex
structure $\theta'$, which differs stably from the first by a unique
difference element $\varDelta(\theta',\theta)$ in $\KO^{-2}(N)$. As
above, its image under $\chi$ is represented by the virtual bundle
$\theta'-\theta$ in $K^0(N)$. Whether or not $\chi$ is monic,
$\varDelta(\theta',\theta)$ is constructed by expressing the trivial
bundle $\bR^{2L}$ as $\theta\oplus\theta^\perp$ for suitably large
$L$, then taking the complex structure induced by $g'$ on $\theta$ and
by the Hermitian complement of $g$ on $\theta^\perp$. We are
particularly interested in this situation when $g'$ is obtained from
$g$ by complex conjugation; the difference element may then be
described as follows.
\begin{lem}\label{diffconj}
For any complex vector bundle $\xi$ over $N$, the difference element
$\varDelta(\overline{\xi},\xi)$ is given by $r(z(\overline{\xi}-1))$
in $\KO^{-2}(N)$.
\end{lem}
\begin{proof}
It is sufficient to consider the universal bundle $\upsilon$ over a
complex Grassmannian of the form $U(W\oplus W')/U(W)\times U(W')$,
where $W\oplus W'$ is isomorphic to $\bC^L$ for suitably large $L$.
Both $\varDelta(\overline{\upsilon},\upsilon)$ and
$r(z(\overline{\upsilon}-\bC))$ may be represented by maps into
$\varOmega^2\SO(W\oplus W')$, obtained by adjointing Bott's original
periodicity maps. Details of these are in \cite{ca:pgh}, as are the
techniques for proving that the two maps are homotopic.
\end{proof}

For any Bott tower $(M^k:k\leq n)$, we write $o(a,k)$ (or $o(k)$ when
the list $a$ is understood or irrelevant) for the number of homotopy
classes of $U$-structures which arise from the omniorientations of
$M^k$. Thus $1\leq o(k)\leq 2^{2k}$. Applying Lemma \ref{diffconj} and
the splitting \eqref{gamgamdav} to the $U$-structure $\tau_\bC$ of
\eqref{canomni} identifies the corresponding difference elements as
\begin{equation}\label{omnidiffs}
\begin{split}
\sum_{j=1}^k\delta_j&\varDelta(\gamma_j,\overline{\gamma}_j)
+\sum_{j=1}^k\epsilon_j
\varDelta\big((\overline{\gamma}(a_{j-1})-\overline{\gamma}_j),
(\gamma(a_{j-1})-\gamma_j)\big)\;=\\
&\sum_{j=1}^k(\delta_j+\epsilon_j)r(z^2g_j)-
\sum_{j=1}^k\epsilon_jr\Big(z^2\prod_{i<j}(g_i+1)^{a(i,j)}\Big),
\end{split}
\end{equation}
where $\delta_j$ and $\epsilon_j$ are $0$ or $1$ for all $1\leq j\leq
k$.

When $k=1$, these reduce to $0$, $xd_1$ and $2xd_1$ in
$\KO^{-2}(M^1)$, so that $o(1)=3$. When $k=2$, Corollary \ref{remt}
shows that we obtain the same elements, together with their
translates by
\[
xd_2,\quad x(d_2-a(1,2)d_1),\sands x(2d_2-a(1,2)d_1)
\]
when $a(1,2)$ is even, and
\[
n_{2,1},\quad n_{2,1}-a(1,2)xd_1,\sands 2n_{2,1}-a(1,2)xd_1
\]
when $a(1,2)$ is is odd. So $o(a,2)=9$, $10$, $11$, and $12$, as
$a(1,2)=0$, $\pm1$, $\pm2$, and $|a(1,2)|\geq 3$ respectively.

The calculations increase rapidly in complexity for general values of
$a(i,j)$. Nevertheless, certain families of special cases yield
interesting conclusions.
\begin{thm}\label{omnisscs}
For any Bott tower $(M^k:k\leq n)$, we have that
\[
3^k\leq o(k)\leq 3\cdot 4^{k-1}
\]
for each $1\leq k\leq n$. The maximum is attained by any tower for
which the inequality $|a(k-1,k)|\geq 3$ holds for all $k$, and the
minimum by the tower $((\CPO)^k:k\geq 0)$; the tower of bounded flag
manifolds $(B_k:k\geq 0)$ satisfies
\[
o(k)=\sum_{i=0}^{\lceil k/2\rceil}\binom{k+1}{2i}2^{k-i}.
\]
\end{thm}
\begin{proof}
We proceed by induction on $k$, having resolved the cases $k=1$ and
$2$ above. We assume first that $|a(k-1,k)|\geq 3$ for all $k$, and
that $o(k-1)=3\cdot 4^{k-2}$. For $M^k$, the difference elements
\eqref{omnidiffs} consist of pullbacks from $M^{k-1}$, plus their
translates by the three nonzero elements
\begin{equation}\label{helpo}
(\delta_k+\epsilon_k)r(z^2g_k)\,+\;
\delta_kr\Big(z^2\prod_{j<k}(g_j+1)^{a(j,k)}\Big).
\end{equation}
These map to $(\delta_k+\epsilon_k)(g_k-\overline{g}_k)+
\delta_k\big(\prod_{j<k}(g_j+1)^{a(j,k)}-
\prod_{j<k}(\overline{g}_j+1)^{a(j,k)}\big)$ under complexification,
where $-\delta_ka(k-1,k)(g_{k-1}-\overline{g}_{k-1})$ is the only term
involving $g_{k-1}$. It follows that no such translates can result in
coincident difference elements when $|a(k-1,k)|\geq 3$, and the
initial induction is complete.

The tower $((\CPO)^k:k\geq 0)$, on the other hand, has $a(i,j)=0$ for
all values of $i$ and $j$, and is totally even. The translation
elements \eqref{helpo} then reduce to $(\delta_k+\epsilon_k)xd_k$,
creating one coincidence for each element pulled back from
$(\CPO)^{k-1}$; this maximises the possible coincidences, and leads to
$o(k)=3o(k-1)$. So $o(k)=3^k$, represented by the difference elements
$\sum_{j=1}^k\omega_jxd_j$, where $\omega_j=0$, $1$, or $2$ for each
$j$.

The tower $(B_k:k\geq 0)$ has $a(j-1,j)=1$ for all $j<k$, and
$a(i,j)=0$ otherwise. Being terminally odd, we may follow Theorem
\ref{komintwo}, and work with the complexifications
\[
g_k-\overline{g}_k,\;
-(g_{k-1}-\overline{g}_{k-1})+(g_k-\overline{g}_k),\; \text{and
$-(g_{k-1}-\overline{g}_{k-1})+2(g_k-\overline{g}_k)$},
\]
of the translation elements \eqref{helpo}. These yield two
coincidences for each element of the $(k-2)$th stage. In other words,
$o(k)$ satisfies the difference equation $o(k)=4o(k-1)-2o(k-2)$ for
each $k\geq 2$. Using the initial conditions provided by $k=1$ and
$2$, we may then apply standard techniques \cite{kepe:de} to deduce
the required formula. The same arguments work when $a(j-1,j)=-1$ and
$a(i,j)=0$ for $i\neq j-1$.
\end{proof}

We emphasise that these results depend on our initial choice of
orientation for $M^k$, as do Theorems \ref{omniacs} and \ref{omnibds}
below. 

It transpires that the $U$-structure $\tau'$ of \eqref{stanmk} is
amongst those induced by an omniorientation, whose difference element
satisfies $\delta_j=1$ and $\epsilon_j=0$ in \eqref{omnidiffs}, for all
$1\leq j\leq k$.
\begin{thm}\label{thrstrs}
For any Bott tower $(M^k:k\leq n)$, the difference element
$\varDelta(\tau',\tau_\bC)$ is given by $\sum_{j=1}^kr(z^2g_j)$ in
$\KO^{-2}(M^k)$.
\end{thm}
\begin{proof}
We proceed by induction on $k$, choosing $k=0$ as the base case because
the elements in question are both zero.

So we assume that the result is true For $M^{k-1}$, and consider the
construction of $M^k$. We observe that $\tau'$
and $\tau_\bC$ both arise by pulling back the corresponding
$U$-structures on $M^{k-1}$, and adding the bundle of tangents along
the fibres. By induction, the structures on $M^{k-1}$ differ by
$\sum_{j=1}^{k-1}r(z^2g_j)$. Moreover, the tangents along the fibres
pull back from the corresponding bundles along the fibres of the
universal example over $\CPI$. In this case, $\KO^*(\CPI)$ is torsion
free, so that $\chi$ is monic and we may work in $K^0(\CPI)$. The
relevant difference element is therefore $r(z^2u)$, and pulls back to
$r(z^2g_k)$ over $M^k$. Adding the results yields the required formula.
\end{proof}

The structure $\tau_\bC$ is the stabilisation of an almost complex
structure, and we would like to estimate how many others that are
induced by an omniorientation share this property. We recall from
Section \ref{boto} our observation that the Euler characteristic
$e(M^k)$ is $2^k$.

According to Thomas \cite{th:csr}, the structures we seek are
precisely those whose $k$th Chern class coincides with $e(M^k)$, and
therefore with $c_k(\tau_\bC)$. We may compute the latter by combining
\eqref{gamgamdav} with \eqref{canomni} and writing the total Chern
class $c(\tau_\bC)$ as
\[
(1-2x_1)\prod_{j=2}^k(1+a(1,j)x_1+\dots+a(j-1,j)x_{j-1}-2x_j).
\]
We deduce that $c_k(\tau_\bC)=(-2)^kx_1\dots x_k$. This confirms the
value of $e(M^k)$, and shows that the orientation class defined by the
complex structure on the projective form of $M^k$ is the dual of
$(-1)^kx_1\dots x_k$ in $H_{2k}(M^k:\bZ)$.
\begin{thm}\label{omniacs}
For any Bott tower $(M^k:k\leq n)$, the omniorientations induce
$2^{k-1}$ distinct almost complex structures on $M^k$, for each $1\leq
k\leq n$.
\end{thm}
\begin{proof}
We may build up the total Chern class of every $U$-structure on $M^k$
by analogy with the proof of Theorem \ref{omnisscs}; when $k=1$ we
obtain $1-2x_1$, $1+2x_1$ or $1$. Only the first of these has the
required $c_1$, confirming the result for $k=1$.

To obtain the $k$th stage, we multiply the $(k-1)$th stage by one of
the four possible factors
\begin{equation}\label{helpth}
\begin{split}
1\pm\big(a(1,k)x_1+\dots&+a(k-1,k)x_{k-1}\big)\quad\text{or}\\
&1\pm\big(a(1,k)x_1+\dots+a(k-1,k)x_{k-1}-2x_k\big).
\end{split}
\end{equation}
The only way in which the monomial $x_1\dots x_k$ (or any of its
equivalent forms such as $x_k^k$) can occur in the final product is by
selecting one of the latter two factors at this, and every previous,
stage. There are $2^k$ such possibilities in all, distributed equally
between $\pm2^kx_1\dots x_k$. 

It remains only to prove that there are no repetitions amongst the
$2^{k-1}$ products with sign $(-1)^k$. In fact all $2^k$ structures
have distinct $c_1$, as a simple computation shows.
\end{proof}

The relevance of bounded flag manifolds to complex cobordism theory
was first highlighted in \cite{ra:ocb}. Somewhat surprisingly, the
most important $U$-structure from this point of view is $\tau'$, which
bounds. We would therefore like to know how many bounding
$U$-structures arise from the omniorientions of $M^k$. We denote this
number by $b(k)$, and conclude with a brief analysis of its possible
values. 
\begin{thm}\label{omnibds}
For any Bott tower $(M^k:k\leq n)$, we have that
\[
3^{k-1}\leq b(k)\leq 3\cdot 4^{k-1}-2^k
\]
for each $1\leq k\leq n$. The towers $((\CPO)^k:k\geq 0)$ and
$(B_k:k\geq 1)$ satisfy 
\[
b(k)=3^k-2^k\sands
b(k)=\sum_{i=0}^{\lceil k/2\rceil}\binom{k}{2i-1}2^{k-i}
\]
respectively.
\end{thm}
\begin{proof}
The lower bound arises from Theorem \ref{omnisscs} by applying
Szczarba's construction \cite{sz:tbf} to deduce that every
$U$-structure on $M^{k-1}$ lifts to a bounding $U$-structure on
$M^k$. The upper bound arises from the fact that the $k$th Chern
number $c_k[M^k]$ of every bounding $U$-structure is zero. Applying
\eqref{helpth} shows that $c_k[M^k]\neq 0$ for precisely $2^k$
distinct $U$-structures, and the inequality $b(k)\leq 3\cdot
4^{k-1}-2^k$ then follows from Theorem \ref{omnisscs}.

The $3^k$ distinct $U$-structures on $(\CPO)^k$ arise by choosing
one of the three possible structures for each factor $\CPO$; one
bounds, the other two do not. A structure on the product bounds 
precisely when one or more of these $k$ choices bound, yielding
$b(k)=3^k-2^k$. For $B_k$, we note from the proof of Theorem 
\ref{omnisscs} that 
\[
b(k)=2o(k-1)-2o(k-2),
\]
so $b(k)$ satisfies $b(k)=4b(k-1)-2b(k-2)$. But there are no bounding
$U$-structures on a point, and only one on $M^1$; so $b(0)=0$ and
$b(1)=1$. Solving the difference equation gives the required formula.
\end{proof}

Many interesting questions remain to be answered about the r\^ole of
Bott towers in complex cobordism theory. We hope to return to these
in future.

%
%
%
%
%
%
%
%
%


\end{document}